\documentclass[12pt]{amsart}
\usepackage{amsfonts}
\usepackage{amsmath}
\usepackage{amssymb}
\usepackage{graphicx,}
\usepackage{enumerate}
\usepackage{multicol}
\usepackage{mathrsfs}
\usepackage[usenames,dvipsnames]{pstricks}
\usepackage{epsfig}
\usepackage{pst-grad} % For gradients
\usepackage{pst-plot} % For axes
\usepackage[hmargin=3cm,vmargin=2.5cm]{geometry}
\usepackage{rotating}
\usepackage[utf8]{inputenc}
\usepackage{rotating}
\usepackage{pstricks, pst-node, pst-text, pst-3d}
\usepackage{amsbsy}
\usepackage{bm}
\usepackage{ mathrsfs}
\usepackage{graphicx}
\usepackage{latexsym}
\usepackage{lscape}
\usepackage{tikz}

\newtheorem{definition}{Definition}

\DeclareMathOperator{\n}{\mathfrak{n}}

\DeclareMathOperator{\z}{\mathfrak{z}}
\DeclareMathOperator{\vv}{\mathfrak{v}}
\DeclareMathOperator{\End}{End}

\DeclareMathOperator{\Id}{Id}
\DeclareMathOperator{\Cl}{\rm{Cl}}

\DeclareMathOperator{\la}{\langle}
\DeclareMathOperator{\ra}{\rangle}

\title{Structure constants of pseudo $H$-type algebras in some integral bases}

\author{K.~Furutani, I.~Markina}

\thanks{The first author has been partially
supported by 
the Grant-in-aid for
Scientific Research (C) No. 23540251, 
{\it Japan Society for the Promotion of Science}. The second author has been partially supported by the grants of the 
Norwegian Research Council \#239033/F20.}
 


\address{K.~Furutani:  Department of Mathematics, Faculty of Science 
and Technology, Science University of Tokyo, 2641 Yamazaki, Noda, Chiba (278-8510), Japan}
\email{furutani\_kenro@ma.noda.tus.ac.jp}

\address{I.~Markina: Department of Mathematics, University of Bergen, P.O.~Box 7803,
Bergen N-5020, Norway}
\email{irina.markina@uib.no}

\begin{document}
\maketitle

\begin{abstract}
We present the structural constants of low dimensional pseudo $H$-type algebras.
\end{abstract}

%%%%%%%%%%%%%%%%%%%%%%%%%%%%%%%%%%%%%%%%%%%%%%%%

\section{Definition of pseudo $H$-type algebra}\label{sec:PLA}

%%%%%%%%%%%%%%%%%%%%%%%%%%%%%%%%%%%%%%%%%%%%%%%%

Let $\n=\z\oplus\vv$ be a nilpotent graded Lie algebra, $\la.\,,.\ra_{r,s}$ a non-degenerate symmetric bi-linear form on $\z$ of index $(r,s)$ and $\la.\,,.\ra$ a non-degenerate symmetric bi-linear form on
$\vv$ that is of index $(l,l)$ in the case $s>0$ and is of index $(n,0)$ if $s=0$. Define the map $J\colon\z\to\End(\vv)$ by
$$
\la z,[x,y]\ra_{r,s}=\la J_zx,y\ra,\quad z\in\z,\ \ x,y\in\vv.
$$
Then the map $J_z$ is skew symmetric with respect to the bi-linear form $\la.\,,.\ra$ in the following sense
\begin{equation}\label{eq:cl_rep}
\la J_zx,y\ra+\la x,J_zy\ra=0.
\end{equation}
\begin{definition}
A Lie algebra $\mathfrak n=(\n,[.\,,.],\la.\,,.\ra_{r,s}+\la.\,,.\ra)$ is called a pseudo $H$-type Lie algebra if
$$
J_z^2+\la z\,,z\ra\Id_{\vv}=0\quad\text{for all}\quad z\in\z.
$$
\end{definition}
The pseudo $H$-type Lie algebras are in one to one correspondence with the Clifford $\Cl(\z,\la.\,,.\ra_{r,s})$-module structures on $\vv$ admitting a symmetric bi-linear form $\la.\,,.\ra$, satisfying~\eqref{eq:cl_rep}, see~\cite{AFM,Ciatti,GKM}. We use the isomorphism of scalar vector spaces $(\z,\la.\,,.\ra_{r,s})\cong\mathbb R^{r,s}$, $(\vv,\la.\,,.\ra)\cong\mathbb R^{l,l}$ or $(\vv,\la.\,,.\ra)\cong\mathbb R^{n,0}$, and the Clifford algebras $\Cl(\z,\la.\,,.\ra_{r,s})\cong\Cl_{r,s}=\Cl(\mathbb R^{r,s})$. We also write $\n_{r,s}$ for pseudo $H$-type Lie algebras related to the $\Cl_{r,s}$-module of minimal possible dimension (minimal admissible module).

We denote by $z_1,\ldots,z_r,z_{r+1},\ldots,z_{r+s}$ an orthonormal basis for $\mathbb R^{r,s}$ such that 
$$
\la z_i,z_i\ra_{r,s}=
\begin{cases}
1,\quad&\text{if}\quad i=1,\ldots,r,
\\
-1&\text{if}\quad i=r+1,\ldots,r+s.
\end{cases}
$$
We choose the initial vector $v$ on each Clifford module $\vv$ such that $\la v,v\ra=1$. More about the construction of the bases for $\n_{r,s}$ and the isomorphism properties see~\cite{AFM,FM14,FM2}.

%%%%%%%%%%%%%%%%%%%%%%%%%%%%%%%%%%%%%%%%%%%%%%%%

\section{Bases and structure constants for pseudo $H$-type Lie algebras with $r+s=1$}

%%%%%%%%%%%%%%%%%%%%%%%%%%%%%%%%%%%%%%%%%%%%%%%%

The pseudo $H$-type Lie algebras $\n_{1,0}$ and $\n_{0,1}$ are isomorphic~\cite{AFM,FM2}. The minimal admissible module $\vv$ is $2$-dimensional. Let $z_1\in\mathbb R^{1,0}$, $\la z_1,z_1\ra_{1,0}=1$, $J_{z_1}^2=-\Id_{\vv}$ and $v\in\vv$ with $\la v,v\ra=1$. Then the integral basis of $\n_{1,0}$ is
$$
z_1,\quad v_1=v,\quad v_2=J_{z_1}v.
$$
Commutators are given in Table~1
{\small
\begin{table}[h]
\caption{Commutation relations for $\n_{1,0}$ and $\n_{0,1}$}
\centering
\begin{tabular}{| c | c | c |} 
\hline
 $[r, c]$  & $v_1$ & $v_2$ \\
\hline
$ v_1$ & $0$ & $z_1$  \\
\hline
$v_2$ & $-z_1$ & $0$ \\
\hline
\end{tabular}\label{10}
\end{table}   }

%%%%%%%%%%%%%%%%%%%%%%%%%%%%%%%%%%%%%%%%%%%%%%%%

\section{Bases and structure constants for pseudo $H$-type Lie algebras with $r+s=2$}

%%%%%%%%%%%%%%%%%%%%%%%%%%%%%%%%%%%%%%%%%%%%%%%%

%%%%%%%%%%%%%%%%%%%%%%%%%%%%%%%%%%%%%%%%%%%%%%%%

\subsection{ $H$-type Lie algebra $\n_{2,0}$} 

%%%%%%%%%%%%%%%%%%%%%%%%%%%%%%%%%%%%%%%%%%%%%%%%

The minimal admissible module $\vv$ is $4$-dimensional. The basis for $\mathbb R^{2,0}$ is $z_1,z_2$ and $J_{z_i}^2=-\Id_{\vv}$, $i=1,2$. The basis for $\vv$ is
$$
v_1=v,\quad v_2=J_{z_2}J_{z_1}v,\quad v_3=J_{z_1}v, \quad v_4=J_{z_2}v.
$$
{ \small
\begin{table}[h]
\caption{Commutation relations for $\n_{2,0}$ and $\n_{0,2}$}
\centering
\begin{tabular}{| c | c | c | c | c |} 
\hline
 $[r, c]$  & $v_1$ & $v_2$ & $v_3$ & $v_4$ \\
\hline
$v_1$ & $0$ & $0$ & $z_1$ & $z_2$ \\
\hline
$v_2$ & $0$ & $0$ & $-z_2$ & $z_1$ \\
\hline
$v_3$ & $-z_1$ & $z_2$ & $0$ & $0$ \\
\hline
$v_4$ & $-z_2$ & $-z_1$ & $0$ & $0$ \\
\hline
\end{tabular}\label{02}
\end{table} }

%%%%%%%%%%%%%%%%%%%%%%%%%%%%%%%%%%%%%%%%%%%%%%%%

\subsection{ $H$-type Lie algebra $\n_{1,1}$} 

%%%%%%%%%%%%%%%%%%%%%%%%%%%%%%%%%%%%%%%%%%%%%%%%

The minimal admissible module $\vv$ is $4$-dimensional. The basis for $\mathbb R^{1,1}$ is $z_1,z_2$ and $J_{z_1}^2=-\Id_{\vv},\ \ J_{z_2}^2=\Id_{\vv}$. The basis for $\vv$ is
$$
v_1=v,\quad v_2=J_{z_1}J_{z_2}v,\quad v_3=J_{z_1}v, \quad v_4=J_{z_2}v.
$$
{ \small
\begin{table}[h]
\caption{Commutation relations for $\n_{1,1}$}
\centering
\begin{tabular}{| c | c | c | c | c |} 
\hline
 $[r, c]$  & $v_1$ & $v_2$ & $v_3$ & $v_4$ \\
\hline
$v_1$ & $0$ & $0$ & $z_1$ & $z_2$ \\
\hline
$v_2$ & $0$ & $0$ & $z_2$ & $z_1$ \\
\hline
$v_3$ & $-z_1$ & $-z_2$ & $0$ & $0$ \\
\hline
$v_4$ & $-z_2$ & $-z_1$ & $0$ & $0$ \\
\hline
\end{tabular}\label{11}
\end{table} }

%%%%%%%%%%%%%%%%%%%%%%%%%%%%%%%%%%%%%%%%%%%%%%%%

\subsection{ $H$-type Lie algebra $\n_{0,2}$} 

%%%%%%%%%%%%%%%%%%%%%%%%%%%%%%%%%%%%%%%%%%%%%%%%

The Lie algebra $\n_{0,2}$ is isomorphic to the Lie algebra $\n_{2,0}$ and the structural constants are given in Table~2.

%%%%%%%%%%%%%%%%%%%%%%%%%%%%%%%%%%%%%%%%%%%%%%%%

\section{Bases and structure constants for pseudo $H$-type Lie algebras with $r+s=3$}

%%%%%%%%%%%%%%%%%%%%%%%%%%%%%%%%%%%%%%%%%%%%%%%%

%%%%%%%%%%%%%%%%%%%%%%%%%%%%%%%%%%%%%%%%%%%%%%%%

\subsection{ $H$-type Lie algebra $\n_{3,0}$} 

%%%%%%%%%%%%%%%%%%%%%%%%%%%%%%%%%%%%%%%%%%%%%%%%

The minimal admissible module $\vv$ is $4$-dimensional. Basis of $\mathbb R^{3,0}$ is $z_1,z_2,z_3$ and we choose an initial vector $v\in\vv$, $\la v,v\ra=1$ satisfying 
$
J_{z_1}J_{z_2}J_{z_3}v=v
$.
Then the basis of $\vv$ is the following
$$
v_1=v,\quad v_2=J_{z_1}v,\quad v_3=J_{z_2}v, \quad v_4=J_{z_3}v.
$$
{ \small
\begin{table}[h]
\caption{Commutation relations for $\n_{3,0}$}
\centering
\begin{tabular}{| c | c | c | c | c |} 
\hline
 $[r, c]$  & $v_1$ & $v_2$ & $v_3$ & $v_4$ \\
\hline
$v_1$ & $0$ & $z_1$ & $z_2$ & $z_3$ \\
\hline
$v_2$ & $-z_1$ & $0$ & $-z_3$ & $z_2$ \\
\hline
$v_3$ & $-z_2$ & $z_3$ & $0$ & $-z_1$ \\
\hline
$v_4$ & $-z_3$ & $-z_2$ & $z_1$ & $0$ \\
\hline
\end{tabular}\label{30}
\end{table} }

%%%%%%%%%%%%%%%%%%%%%%%%%%%%%%%%%%%%%%%%%%%%%%%%

\subsection{ $H$-type Lie algebra $\n_{2,1}$} 

%%%%%%%%%%%%%%%%%%%%%%%%%%%%%%%%%%%%%%%%%%%%%%%%

The minimal admissible module $\vv$ is $8$-dimensional. Basis of $\mathbb R^{2,1}$ is $z_1,z_2,z_3$ and $J^2_{z_i}=-\Id_{\vv}$, $i=1,2$, $J^2_{z_3}=\Id_{\vv}$. We choose an initial vector $v\in\vv$ such that $\la v, J_{z_1}J_{z_2}J_{z_3}v\ra=0$ and $\la v,v\ra=1$. Then the basis of $\vv$ is the following
$$
\begin{array}{lllllll}
v_1=v, \quad &v_2=J_{z_2}J_{z_1}v,\quad &v_3=J_{z_1}J_{z_3}v, \quad &v_4=J_{z_2}J_{z_3}v,
\\
v_5=J_{z_1}v,\quad &v_6=J_{z_2}v,\quad &v_7=J_{z_3}v, \quad &v_8=J_{z_1}J_{z_2}J_{z_3}v.
\end{array}
$$
{ \small
\begin{table}[h]
\caption{Commutation relations for $\n_{2,1}$}
\centering
\begin{tabular}{| c | c | c | c | c | c | c | c | c | c |} 
\hline
 $[r, c]$  & $v_1$ & $v_2$ & $v_3$ & $v_4$ & $v_5$ & $v_6$ & $v_7$ & $v_8$
 \\
\hline
$v_1$ & $0$      & $0$      & $0$      & $0$ & $z_1$ & $z_2$ & $z_3$ & $0$\\
\hline
$v_2$ & $0$      & $0$      & $0$      & $0$ & $-z_2$ & $z_1$ & $0$ & $-z_3$ \\
\hline
$v_3$ & $0$      & $0$      & $0$      & $0$ & $z_3$ & $0$ & $z_1$ & $z_2$ \\
\hline
$v_4$ & $0$      & $0$      & $0$      & $0$ & $0$ & $z_3$ & $z_2$ & $-z_1$\\
\hline
$v_5$ & $-z_1$ & $z_2$ & $-z_3$ & $0$ & $0$ & $0$ & $0$ & $0$\\
\hline
$v_6$ & $-z_2$ & $-z_1$  & $0$      & $-z_3$ & $0$ & $0$ & $0$ & $0$\\
\hline
$v_7$ & $-z_3$ & $0$      & $-z_1$ & $-z_2$ & $0$ & $0$ & $0$ & $0$\\
\hline
$v_8$ & $0$      & $z_3$  & $-z_2$ & $z_1$ & $0$ & $0$ & $0$ & $0$\\
\hline
\end{tabular}\label{21-1}
\end{table} }

%%%%%%%%%%%%%%%%%%%%%%%%%%%%%%%%%%%%%%%%%%%%%%%%

\subsection{ $H$-type Lie algebra $\n_{1,2}$} 

%%%%%%%%%%%%%%%%%%%%%%%%%%%%%%%%%%%%%%%%%%%%%%%%

The minimal admissible module $\vv$ is $4$-dimensional. Basis of $\mathbb R^{1,2}$ is $z_1,z_2,z_3$ and $J^2_{z_1}=-\Id_{\vv}$, $J^2_{z_i}=\Id_{\vv}$, $i=2,3$.  We choose an initial vector $v\in\vv$, $\la v,v\ra=1$, in order to satisfy 
$
J_{z_1}J_{z_2}J_{z_3}v=v$.
Then the basis of $\vv$ is the following 
$$
v_1=v,\quad v_2=J_{z_1}v,\quad v_3=J_{z_2}v, \quad v_4=J_{z_3}v.
$$
{ \small
\begin{table}[h]
\caption{Commutation relations for $\n_{1,2}$}
\centering
\begin{tabular}{| c | c | c | c | c |} 
\hline
 $[r, c]$  & $v_1$ & $v_2$ & $v_3$ & $v_4$ \\
\hline
$v_1$ & $0$ & $z_1$ & $z_2$ & $z_3$ \\
\hline
$v_2$ & $-z_1$ & $0$ & $z_3$ & $-z_2$ \\
\hline
$v_3$ & $-z_2$ & $-z_3$ & $0$ & $-z_1$ \\
\hline
$v_4$ & $-z_3$ & $z_2$ & $z_1$ & $0$ \\
\hline
\end{tabular}\label{12}
\end{table} }

%%%%%%%%%%%%%%%%%%%%%%%%%%%%%%%%%%%%%%%%%%%%%%%%

\subsection{ $H$-type Lie algebra $\n_{0,3}$} 

%%%%%%%%%%%%%%%%%%%%%%%%%%%%%%%%%%%%%%%%%%%%%%%%

The minimal admissible module $\vv$ is $8$-dimensional. Basis of $\mathbb R^{0,3}$ is $z_1,z_2,z_3$ and $J^2_{z_i}=\Id_{\vv}$, $i=1,2,3$. We choose an initial vector $v\in\vv$, $\la v,v\ra=1$, such that $\la v, J_{z_1}J_{z_2}J_{z_3}v\ra=0$. Then the basis of $\vv$ is the following
$$
\begin{array}{lllllll}
v_1=v, \quad &v_2=J_{z_1}J_{z_2}v,\quad &v_3=J_{z_1}J_{z_3}v, \quad &v_4=J_{z_2}J_{z_3}v,
\\
v_5=J_{z_1}v,\quad &v_6=J_{z_2}v,\quad &v_7=J_{z_3}v, \quad &v_8=J_{z_1}J_{z_2}J_{z_3}v
\end{array}
$$
{ \small
\begin{table}[h]
\caption{Commutation relations for $\n_{0,3}$}
\centering
\begin{tabular}{| c | c | c | c | c | c | c | c | c | c |} 
\hline
 $[r, c]$  & $v_1$ & $v_2$ & $v_3$ & $v_4$ & $v_5$ & $v_6$ & $v_7$ & $v_8$
 \\
\hline
$v_1$ & $0$      & $0$      & $0$      & $0$ & $z_1$ & $z_2$ & $z_3$ & $0$\\
\hline
$v_2$ & $0$      & $0$      & $0$      & $0$ & $-z_2$ & $z_1$ & $0$ & $z_3$ \\
\hline
$v_3$ & $0$      & $0$      & $0$      & $0$ & $-z_3$ & $0$ & $z_1$ & $-z_2$ \\
\hline
$v_4$ & $0$      & $0$      & $0$      & $0$ & $0$ & $-z_3$ & $z_2$ & $z_1$\\
\hline
$v_5$ & $-z_1$ & $z_2$ & $z_3$ & $0$ & $0$ & $0$ & $0$ & $0$\\
\hline
$v_6$ & $-z_2$ & $-z_1$  & $0$      & $z_3$ & $0$ & $0$ & $0$ & $0$\\
\hline
$v_7$ & $-z_3$ & $0$      & $-z_1$ & $-z_2$ & $0$ & $0$ & $0$ & $0$\\
\hline
$v_8$ & $0$      & $-z_3$  & $z_2$ & $-z_1$ & $0$ & $0$ & $0$ & $0$\\
\hline
\end{tabular}\label{03}
\end{table} }

%%%%%%%%%%%%%%%%%%%%%%%%%%%%%%%%%%%%%%%%%%%%%%%%

\section{Bases and structure constants for pseudo $H$-type Lie algebras with $r+s=4$}

%%%%%%%%%%%%%%%%%%%%%%%%%%%%%%%%%%%%%%%%%%%%%%%%

%%%%%%%%%%%%%%%%%%%%%%%%%%%%%%%%%%%%%%%%%%%%%%%%

\subsection{ $H$-type Lie algebra $\n_{4,0}$} 

%%%%%%%%%%%%%%%%%%%%%%%%%%%%%%%%%%%%%%%%%%%%%%%%

The minimal admissible module $\vv$ is $8$-dimensional. Basis of $\mathbb R^{4,0}$ is $z_1,\ldots,z_4$ and $J_{z_i}^2=-\Id_{\vv}$, $i=1,\ldots,4$. We choose an initial vector $v\in\vv$, $\la v,v\ra=1$, in order to satisfy 
$
J_{z_1}J_{z_2}J_{z_3}J_{z_4}v=v$.
Then the basis of $\vv$ is the following
$$
\begin{array}{llllllll}
&v_1=v,\quad &v_2=J_{z_1}J_{z_2}v,\quad &v_3=J_{z_1}J_{z_3}v, \quad &v_4=J_{z_1}J_{z_4}v,
\\
&v_5=J_{z_1}v,\quad &v_6=J_{z_2}v,\quad &v_7=J_{z_3}v, \quad &v_8=J_{z_4}v.
\end{array}
$$
{ \small
\begin{table}[h]
\caption{Commutation relations for $\n_{4,0}$}
\centering
\begin{tabular}{| c | c | c | c | c | c | c | c | c | c |} 
\hline
 $[r, c]$  & $v_1$ & $v_2$ & $v_3$ & $v_4$ & $v_5$ & $v_6$ & $v_7$ & $v_8$
 \\
\hline
$v_1$ & $0$      & $0$      & $0$      & $0$ & $z_1$ & $z_2$ & $z_3$ & $z_4$\\
\hline
$v_2$ & $0$      & $0$      & $0$      & $0$ & $z_2$ & $-z_1$ & $-z_4$ & $z_3$ \\
\hline
$v_3$ & $0$      & $0$      & $0$      & $0$ & $z_3$ & $z_4$ & $-z_1$ & $-z_2$ \\
\hline
$v_4$ & $0$      & $0$      & $0$      & $0$ & $z_4$ & $-z_3$ & $z_2$ & $-z_1$\\
\hline
$v_5$ & $-z_1$ & $-z_2$ & $-z_3$ & $-z_4$ & $0$ & $0$ & $0$ & $0$\\
\hline
$v_6$ & $-z_2$ & $z_1$  & $-z_4$ & $z_3$ & $0$ & $0$ & $0$ & $0$\\
\hline
$v_7$ & $-z_3$ & $z_4$  & $z_1$ & $-z_2$ & $0$ & $0$ & $0$ & $0$\\
\hline
$v_8$ & $-z_4$ & $-z_3$  & $z_2$ & $z_1$ & $0$ & $0$ & $0$ & $0$\\
\hline
\end{tabular}\label{40}
\end{table} }

%%%%%%%%%%%%%%%%%%%%%%%%%%%%%%%%%%%%%%%%%%%%%%%%

\subsection{ $H$-type Lie algebra $\n_{3,1}$} 

%%%%%%%%%%%%%%%%%%%%%%%%%%%%%%%%%%%%%%%%%%%%%%%%

The minimal admissible module $\vv$ is $8$-dimensional. Basis of $\mathbb R^{3,1}$ is $z_1,\ldots,z_4$ and $J_{z_i}^2=-\Id_{\vv}$, $i=1,2,3$, $J_{z_4}^2=\Id_{\vv}$. We choose an initial vector $v\in\vv$, $\la v,v\ra=1$, such that  
$
J_{z_1}J_{z_2}J_{z_3}v=v$.
Then the basis of $\vv$ is the following
$$
\begin{array}{llllllll}
&v_1=v,\quad &v_2=J_{z_1}v,\quad &v_3=J_{z_2}v, \quad &v_4=J_{z_3}v,
\\
&v_5=J_{z_4}v,\quad &v_6=J_{z_4}J_{z_1}v,\quad &v_7=J_{z_4}J_{z_2}v, \quad &v_8=J_{z_4}J_{z_3}v.
\end{array}
$$
{ \small
\begin{table}[h]
\caption{Commutation relations for $\n_{3,1}$}
\centering
\begin{tabular}{| c | c | c | c | c | c | c | c | c | c |} 
\hline
 $[r, c]$  & $v_1$ & $v_2$ & $v_3$ & $v_4$ & $v_5$ & $v_6$ & $v_7$ & $v_8$
 \\
\hline
$v_1$ & $0$ & $z_1$ & $z_2$ & $z_3$ & $z_4$ & $0$ & $0$ & $0$\\
\hline
$v_2$ & $-z_1$ & $0$ & $-z_3$ & $z_2$ & $0$ & $z_4$ & $0$ & $0$ \\
\hline
$v_3$ & $-z_2$ & $z_3$ & $0$ & $-z_1$ & $0$ & $0$ & $z_4$ & $0$ \\
\hline
$v_4$ & $-z_3$ & $-z_2$ & $z_1$ & $0$ & $0$ & $0$ & $0$ & $z_4$\\
\hline
$v_5$ & $-z_4$ & $0$ & $0$ & $0$ & $0$ & $z_1$ & $z_2$ & $z_3$\\
\hline
$v_6$ & $0$ & $-z_4$  & $0$ & $0$ & $-z_1$ & $0$ & $-z_3$ & $z_2$\\
\hline
$v_7$ & $0$ & $0$  & $-z_4$ & $0$ & $-z_2$ & $z_3$ & $0$ & $-z_1$\\
\hline
$v_8$ & $0$ & $0$  & $0$ & $-z_4$ & $-z_3$ & $-z_2$ & $z_1$ & $0$\\
\hline
\end{tabular}\label{31}
\end{table} }

%%%%%%%%%%%%%%%%%%%%%%%%%%%%%%%%%%%%%%%%%%%%%%%%

\subsection{ $H$-type Lie algebra $\n_{2,2}$} 

%%%%%%%%%%%%%%%%%%%%%%%%%%%%%%%%%%%%%%%%%%%%%%%%

The minimal admissible module $\vv$ is $8$-dimensional. Basis of $\mathbb R^{2,2}$ is $z_1,\ldots,z_4$ and $J_{z_i}^2=-\Id_{\vv}$, $i=1,2$, $J_{z_i}^2=\Id_{\vv}$, $i=3,4$. We choose an initial vector $v\in\vv$, $\la v,v\ra=1$ such that  
$
J_{z_1}J_{z_2}J_{z_3}J_{z_4}v=v$.
Then the basis of $\vv$ is the following 
$$
\begin{array}{llllllll}
&v_1=v,\quad &v_2=J_{z_1}J_{z_2}v,\quad &v_3=J_{z_1}J_{z_3}v, \quad &v_4=J_{z_1}J_{z_4}v,
\\
&v_5=J_{z_1}v,\quad &v_6=J_{z_2}v,\quad &v_7=J_{z_3}v, \quad &v_8=J_{z_4}v.
\end{array}
$$
{ \small
\begin{table}[h]
\caption{Commutation relations for $\n_{2,2}$}
\centering
\begin{tabular}{| c | c | c | c | c | c | c | c | c | c |} 
\hline
 $[r, c]$  & $v_1$ & $v_2$ & $v_3$ & $v_4$ & $v_5$ & $v_6$ & $v_7$ & $v_8$
 \\
\hline
$v_1$ & $0$      & $0$      & $0$      & $0$ & $z_1$ & $z_2$ & $z_3$ & $z_4$\\
\hline
$v_2$ & $0$      & $0$      & $0$      & $0$ & $z_2$ & $-z_1$ & $z_4$ & $-z_3$ \\
\hline
$v_3$ & $0$      & $0$      & $0$      & $0$ & $z_3$ & $-z_4$ & $z_1$ & $-z_2$ \\
\hline
$v_4$ & $0$      & $0$      & $0$      & $0$ & $z_4$ & $z_3$ & $z_2$ & $z_1$\\
\hline
$v_5$ & $-z_1$ & $-z_2$ & $-z_3$ & $-z_4$ & $0$ & $0$ & $0$ & $0$\\
\hline
$v_6$ & $-z_2$ & $z_1$  & $z_4$ & $-z_3$ & $0$ & $0$ & $0$ & $0$\\
\hline
$v_7$ & $-z_3$ & $-z_4$  & $-z_1$ & $-z_2$ & $0$ & $0$ & $0$ & $0$\\
\hline
$v_8$ & $-z_4$ & $-z_3$  & $z_2$ & $-z_1$ & $0$ & $0$ & $0$ & $0$\\
\hline
\end{tabular}\label{22}
\end{table} }

%%%%%%%%%%%%%%%%%%%%%%%%%%%%%%%%%%%%%%%%%%%%%%%%

\subsection{ $H$-type Lie algebra $\n_{1,3}$} 

%%%%%%%%%%%%%%%%%%%%%%%%%%%%%%%%%%%%%%%%%%%%%%%%

the minimal admissible module $\vv$ is $8$-dimensional. Basis of $\mathbb R^{1,3}$ is $z_1,\ldots,z_4$ and $J_{z_1}^2=-\Id_{\vv}$, $J_{z_i}^2=\Id_{\vv}$, $i=2,3,4$. We choose an initial vector $v\in\vv$, $\la v,v\ra=1$ satisfying 
$
J_{z_1}J_{z_2}J_{z_3}v=v$.
Then the basis of $\vv$ is the following
$$
\begin{array}{llllllll}
&v_1=v,\quad &v_2=J_{z_1}v,\quad &v_3=J_{z_2}v, \quad &v_4=J_{z_3}v,
\\
&v_5=J_{z_4}v,\quad &v_6=J_{z_4}J_{z_1}v,\quad &v_7=J_{z_4}J_{z_2}v, \quad &v_8=J_{z_4}J_{z_3}v.
\end{array}
$$
{ \small
\begin{table}[h]
\caption{Commutation relations for $\n_{1,3}$}
\centering
\begin{tabular}{| c | c | c | c | c | c | c | c | c | c |} 
\hline
 $[r, c]$  & $v_1$ & $v_2$ & $v_3$ & $v_4$ & $v_5$ & $v_6$ & $v_7$ & $v_8$
 \\
\hline
$v_1$ & $0$ & $z_1$ & $z_2$ & $z_3$ & $z_4$ & $0$ & $0$ & $0$\\
\hline
$v_2$ & $-z_1$ & $0$ & $z_3$ & $-z_2$ & $0$ & $z_4$ & $0$ & $0$ \\
\hline
$v_3$ & $-z_2$ & $-z_3$ & $0$ & $-z_1$ & $0$ & $0$ & $-z_4$ & $0$ \\
\hline
$v_4$ & $-z_3$ & $z_2$ & $z_1$ & $0$ & $0$ & $0$ & $0$ & $-z_4$\\
\hline
$v_5$ & $-z_4$ & $0$ & $0$ & $0$ & $0$ & $z_1$ & $z_2$ & $z_3$\\
\hline
$v_6$ & $0$ & $-z_4$  & $0$ & $0$ & $-z_1$ & $0$ & $z_3$ & $-z_2$\\
\hline
$v_7$ & $0$ & $0$  & $z_4$ & $0$ & $-z_2$ & $-z_3$ & $0$ & $-z_1$\\
\hline
$v_8$ & $0$ & $0$  & $0$ & $z_4$ & $-z_3$ & $z_2$ & $z_1$ & $0$\\
\hline
\end{tabular}\label{13}
\end{table} }

%%%%%%%%%%%%%%%%%%%%%%%%%%%%%%%%%%%%%%%%%%%%%%%%

\subsection{ $H$-type Lie algebra $\n_{0,4}$} 

%%%%%%%%%%%%%%%%%%%%%%%%%%%%%%%%%%%%%%%%%%%%%%%%

The minimal admissible module $\vv$ is $8$-dimensional. Basis of $\mathbb R^{0,4}$ is $z_1,\ldots,z_4$ and $J_{z_i}^2=\Id_{\vv}$, $i=1,\ldots,4$. We choose an initial vector $v\in\vv$, $\la v,v\ra=1$, such that  
$
J_{z_1}J_{z_2}J_{z_3}J_{z_4}v=v$. 
Then the basis of $\vv$ is the following
$$
\begin{array}{llllllll}
&v_1=v,\quad &v_2=J_{z_2}J_{z_1}v,\quad &v_3=J_{z_3}J_{z_1}v, \quad &v_4=J_{z_4}J_{z_1}v,
\\
&v_5=J_{z_1}v,\quad &v_6=J_{z_2}v,\quad &v_7=J_{z_3}v, \quad &v_8=J_{z_4}v.
\end{array}
$$
{ \small
\begin{table}[h]
\caption{Commutation relations for $\n_{0,4}$}
\centering
\begin{tabular}{| c | c | c | c | c | c | c | c | c | c |} 
\hline
 $[r, c]$  & $v_1$ & $v_2$ & $v_3$ & $v_4$ & $v_5$ & $v_6$ & $v_7$ & $v_8$
 \\
\hline
$v_1$ & $0$      & $0$      & $0$      & $0$ & $z_1$ & $z_2$ & $z_3$ & $z_4$\\
\hline
$v_2$ & $0$      & $0$      & $0$      & $0$ & $z_2$ & $-z_1$ & $-z_4$ & $z_3$ \\
\hline
$v_3$ & $0$      & $0$      & $0$      & $0$ & $z_3$ & $z_4$ & $-z_1$ & $-z_2$ \\
\hline
$v_4$ & $0$      & $0$      & $0$      & $0$ & $z_4$ & $-z_3$ & $z_2$ & $-z_1$\\
\hline
$v_5$ & $-z_1$ & $-z_2$ & $-z_3$ & $-z_4$ & $0$ & $0$ & $0$ & $0$\\
\hline
$v_6$ & $-z_2$ & $z_1$  & $-z_4$ & $z_3$ & $0$ & $0$ & $0$ & $0$\\
\hline
$v_7$ & $-z_3$ & $z_4$  & $z_1$ & $-z_2$ & $0$ & $0$ & $0$ & $0$\\
\hline
$v_8$ & $-z_4$ & $-z_3$  & $z_2$ & $z_1$ & $0$ & $0$ & $0$ & $0$\\
\hline
\end{tabular}\label{04}
\end{table} }

We see that Lie algebras $\n_{4,0}$ and $\n_{0,4}$ are isomorphic.

%%%%%%%%%%%%%%%%%%%%%%%%%%%%%%%%%%%%%%%%%%%%%%%%

\section{Bases and structure constants for pseudo $H$-type Lie algebras with $r+s=5$}

%%%%%%%%%%%%%%%%%%%%%%%%%%%%%%%%%%%%%%%%%%%%%%%%

%%%%%%%%%%%%%%%%%%%%%%%%%%%%%%%%%%%%%%%%%%%%%%%%

\subsection{ $H$-type Lie algebra $\n_{5,0}$} 

%%%%%%%%%%%%%%%%%%%%%%%%%%%%%%%%%%%%%%%%%%%%%%%%

The minimal admissible module $\vv$ is $8$-dimensional. Basis of $\mathbb R^{5,0}$ is $z_1,\ldots,z_5$ and $J_{z_i}^2=-\Id_{\vv}$, $i=1,\ldots,5$. We choose an initial vector $v\in\vv$, $\la v,v\ra=1$ to satisfy 
$$
P_1v=J_{z_1}J_{z_2}J_{z_3}J_{z_4}v=v,\quad P_2v=J_{z_1}J_{z_2}J_{z_5}v=v.
$$
Then the basis of $\vv$ is the following
$$
\begin{array}{llllllll}
&v_1=v,\quad &v_2=J_{z_5}v,\quad &v_3=J_{z_1}J_{z_3}v, \quad &v_4=J_{z_1}J_{z_4}v,
\\
&v_5=J_{z_1}v,\quad &v_6=J_{z_2}v,\quad &v_7=J_{z_3}v, \quad &v_8=J_{z_4}v.
\end{array}
$$
{ \small
\begin{table}[h]
\caption{Commutation relations for $\n_{5,0}$}
\centering
\begin{tabular}{| c | c | c | c | c | c | c | c | c | c |} 
\hline
 $[r, c]$  & $v_1$ & $v_2$ & $v_3$ & $v_4$ & $v_5$ & $v_6$ & $v_7$ & $v_8$
 \\
\hline
$v_1$ & $0$ & $z_5$ & $0$ & $0$ & $z_1$ & $z_2$ & $z_3$ & $z_4$\\
\hline
$v_2$ & $-z_5$ & $0$ & $0$ & $0$ & $-z_2$ & $z_1$ & $z_4$ & $-z_3$ \\
\hline
$v_3$ & $0$ & $0$ & $0$ & $-z_5$ & $z_3$ & $z_4$ & $-z_1$ & $-z_2$ \\
\hline
$v_4$ & $0$ & $0$ & $z_5$ & $0$ & $z_4$ & $-z_3$ & $z_2$ & $-z_1$\\
\hline
$v_5$ & $-z_1$ & $z_2$ & $-z_3$ & $-z_4$ & $0$ & $-z_5$ & $0$ & $0$\\
\hline
$v_6$ & $-z_2$ & $-z_1$  & $-z_4$ & $z_3$ & $z_5$ & $0$ & $0$ & $0$\\
\hline
$v_7$ & $-z_3$ & $-z_4$  & $z_1$ & $-z_2$ & $0$ & $0$ & $0$ & $z_5$\\
\hline
$v_8$ & $-z_4$ & $z_3$  & $z_2$ & $z_1$ & $0$ & $0$ & $-z_5$ & $0$\\
\hline
\end{tabular}\label{50}
\end{table} }

%%%%%%%%%%%%%%%%%%%%%%%%%%%%%%%%%%%%%%%%%%%%%%%%

\subsection{ $H$-type Lie algebra $\n_{4,1}$} 

%%%%%%%%%%%%%%%%%%%%%%%%%%%%%%%%%%%%%%%%%%%%%%%%

The minimal admissible module $\vv$ is $16$-dimensional. Basis of $\mathbb R^{4,1}$ is $z_1,\ldots,z_5$ and $J_{z_i}^2=-\Id_{\vv}$, $i=1,\ldots,4$, $J_{z_5}^2=\Id_{\vv}$. We choose an initial vector $v\in\vv$, $\la v,v\ra=1$, satisfying 
$
P_1v=J_{z_1}J_{z_2}J_{z_3}J_{z_4}v=v$.
Then the basis of $\vv$ is the following
$$
\begin{array}{lllllllllllllll}
&v_1=v,\quad &v_2=J_{z_1}J_{z_2}v,\quad &v_3=J_{z_1}J_{z_3}v, \quad &v_4=J_{z_1}J_{z_4}v,
\\
&v_5=J_{z_1}v,\quad &v_6=J_{z_2}v,\quad &v_7=J_{z_3}v, \quad &v_8=J_{z_4}v,
\\
&v_9=J_{z_5}v,\quad &v_{10}=J_{z_5}J_{z_1}J_{z_2}v,\quad &v_{11}=J_{z_5}J_{z_1}J_{z_3}v, \quad &v_{12}=J_{z_5}J_{z_1}J_{z_4}v,
\\
&v_{13}=J_{z_5}J_{z_1}v,\quad &v_{14}=J_{z_5}J_{z_2}v,\quad &v_{15}=J_{z_5}J_{z_3}v, \quad &v_{16}=J_{z_5}J_{z_4}v.
\end{array}
$$
{ \tiny
\begin{table}[h]
\caption{Commutation relations for $\n_{4,1}$}
\centering
\begin{tabular}{| c | c | c | c | c | c | c | c | c | c | c | c | c | c | c | c | c | c | c | c |} 
\hline
 $[r, c]$  & $v_1$ & $v_2$ & $v_3$ & $v_4$ & $v_5$ & $v_6$ & $v_7$ & $v_8$ & $v_9$ & $v_{10}$ & $v_{11}$ & $v_{12}$ & $v_{13}$ & $v_{14}$ & $v_{15}$ & $v_{16}$
 \\
\hline
$v_1$ & $0$ & $0$ & $0$ & $0$ & $z_1$ & $z_2$ & $z_3$ & $z_4$ & $z_5$ & $0$ & $0$ & $0$ & $0$ & $0$ & $0$ & $0$\\
\hline
$v_2$ & $0$ & $0$ & $0$ & $0$ & $z_2$ & $-z_1$ & $-z_4$ & $z_3$ & $0$ & $z_5$ & $0$ & $0$ & $0$ & $0$ & $0$ & $0$\\
\hline
$v_3$ & $0$ & $0$ & $0$ & $0$ & $z_3$ & $z_4$ & $-z_1$ & $-z_2$ & $0$ & $0$ & $z_5$ & $0$ & $0$ & $0$ & $0$ & $0$\\
\hline
$v_4$ & $0$ & $0$ & $0$ & $0$ & $z_4$ & $-z_3$ & $z_2$ & $-z_1$ & $0$ & $0$ & $0$ & $z_5$ & $0$ & $0$ & $0$ & $0$\\
\hline
$v_5$ & $-z_1$ & $-z_2$ & $-z_3$ & $-z_4$ & $0$ & $0$ & $0$ & $0$ & $0$ & $0$ & $0$ & $0$ & $z_5$ & $0$ & $0$ & $0$\\
\hline
$v_6$ & $-z_2$ & $z_1$  & $-z_4$ & $z_3$ & $0$ & $0$ & $0$ & $0$ & $0$ & $0$ & $0$ & $0$ & $0$ & $z_5$ & $0$ & $0$\\
\hline
$v_7$ & $-z_3$ & $z_4$  & $z_1$ & $-z_2$ & $0$ & $0$ & $0$ & $0$ & $0$ & $0$ & $0$ & $0$ & $0$ & $0$ & $z_5$ & $0$\\
\hline
$v_8$ & $-z_4$ & $-z_3$  & $z_2$ & $z_1$ & $0$ & $0$ & $0$ & $0$ & $0$ & $0$ & $0$ & $0$ & $0$ & $0$ & $0$ & $z_5$\\
\hline
$v_9$ & $-z_5$ & $0$  & $0$ & $0$ & $0$ & $0$ & $0$ & $0$ & $0$ & $0$  & $0$ & $0$ & $z_1$ & $z_2$ & $z_3$ & $z_4$\\
\hline
$v_{10}$ & $0$ & $-z_5$  & $0$ & $0$ & $0$ & $0$ & $0$ & $0$ & $0$ & $0$  & $0$ & $0$ & $z_2$ & $-z_1$ & $-z_4$ & $z_3$\\
\hline
$v_{11}$ & $0$ & $0$  & $-z_5$ & $0$ & $0$ & $0$ & $0$ & $0$ & $0$ & $0$  & $0$ & $0$ & $z_3$ & $z_4$ & $-z_1$ & $-z_2$\\
\hline
$v_{12}$ & $0$ & $0$  & $0$ & $-z_5$ & $0$ & $0$ & $0$ & $0$ & $0$ & $0$  & $0$ & $0$ & $z_4$ & $-z_3$ & $z_2$ & $-z_1$ \\
\hline
$v_{13}$ & $0$ & $0$  & $0$ & $0$ & $-z_5$ & $0$ & $0$ & $0$ & $-z_1$ & $-z_2$  & $-z_3$ & $-z_4$ & $0$ & $0$ & $0$ & $0$\\
\hline
$v_{14}$ & $0$ & $0$  & $0$ & $0$ & $0$ & $-z_5$ & $0$ & $0$ & $-z_2$ & $z_1$  & $-z_4$ & $z_3$ & $0$ & $0$ & $0$ & $0$\\
\hline
$v_{15}$ & $0$ & $0$  & $0$ & $0$ & $0$ & $0$ & $-z_5$ & $0$ & $-z_3$ & $z_4$  & $z_1$ & $-z_2$ & $0$ & $0$ & $0$ & $0$\\
\hline
$v_{16}$ & $0$ & $0$  & $0$ & $0$ & $0$ & $0$ & $0$ & $-z_5$ & $-z_4$ & $-z_3$  & $z_2$ & $z_1$ & $0$ & $0$ & $0$ & $0$\\
\hline
\end{tabular}\label{41}
\end{table} }

%%%%%%%%%%%%%%%%%%%%%%%%%%%%%%%%%%%%%%%%%%%%%%%%

\subsection{ $H$-type Lie algebra $\n_{3,2}$} 

%%%%%%%%%%%%%%%%%%%%%%%%%%%%%%%%%%%%%%%%%%%%%%%%

The minimal admissible module $\vv$ is $8$-dimensional. Basis of $\mathbb R^{3,2}$ is $z_1,\ldots,z_5$ and $J_{z_i}^2=-\Id_{\vv}$, $i=1,2,3$, $J_{z_i}^2=\Id_{\vv}$, $i=4,5$. We choose an initial vector $v\in\vv$, $\la v,v\ra=1$, such that  
$$
P_1v=J_{z_2}J_{z_3}J_{z_4}J_{z_5}v=v,\quad P_2v=J_{z_1}J_{z_2}J_{z_3}v=v.
$$
Then the basis of $\vv$ is the following
$$
\begin{array}{llllllll}
&v_1=v,\quad &v_2=J_{z_1}v,\quad &v_3=J_{z_2}v, \quad &v_4=J_{z_3}v,
\\
&v_5=J_{z_4}v,\quad &v_6=J_{z_5}v,\quad &v_7=J_{z_4}J_{z_2}v, \quad &v_8=J_{z_4}J_{z_3}v.
\end{array}
$$
{ \small
\begin{table}[h]
\caption{Commutation relations for $\n_{3,2}$}
\centering
\begin{tabular}{| c | c | c | c | c | c | c | c | c | c |} 
\hline
 $[r, c]$  & $v_1$ & $v_2$ & $v_3$ & $v_4$ & $v_5$ & $v_6$ & $v_7$ & $v_8$
 \\
\hline
$v_1$ & $0$ & $z_1$ & $z_2$ & $z_3$ & $z_4$ & $z_5$ & $0$ & $0$\\
\hline
$v_2$ & $-z_1$ & $0$ & $-z_3$ & $z_2$ & $-z_5$ & $z_4$ & $0$ & $0$ \\
\hline
$v_3$ & $-z_2$ & $z_3$ & $0$ & $-z_1$ & $0$ & $0$ & $z_4$ & $z_5$ \\
\hline
$v_4$ & $-z_3$ & $-z_2$ & $z_1$ & $0$ & $0$ & $0$ & $-z_5$ & $z_4$\\
\hline
$v_5$ & $-z_4$ & $z_5$ & $0$ & $0$ & $0$ & $z_1$ & $z_2$ & $z_3$\\
\hline
$v_6$ & $-z_5$ & $-z_4$  & $0$ & $0$ & $-z_1$ & $0$ & $-z_3$ & $z_2$\\
\hline
$v_7$ & $0$ & $0$  & $-z_4$ & $-z_5$ & $-z_2$ & $z_3$ & $0$ & $-z_1$\\
\hline
$v_8$ & $0$ & $0$  & $z_5$ & $-z_4$ & $-z_3$ & $-z_2$ & $z_1$ & $0$\\
\hline
\end{tabular}\label{32}
\end{table} }

%%%%%%%%%%%%%%%%%%%%%%%%%%%%%%%%%%%%%%%%%%%%%%%%

\subsection{ $H$-type Lie algebra $\n_{2,3}$} 

%%%%%%%%%%%%%%%%%%%%%%%%%%%%%%%%%%%%%%%%%%%%%%%%

The minimal admissible module $\vv$ is $8$-dimensional. Basis of $\mathbb R^{2,3}$ is $z_1,\ldots,z_5$ and $J_{z_i}^2=-\Id_{\vv}$, $i=1,2$, $J_{z_i}^2=\Id_{\vv}$, $i=3,4,5$. We choose an initial vector $v\in\vv$, $\la v,v\ra=1$, such that 
$$
P_1v=J_{z_1}J_{z_2}J_{z_3}J_{z_4}v=v,\quad P_2v=J_{z_1}J_{z_4}J_{z_5}v=v. 
$$
Then the basis of $\vv$ is the following
$$
\begin{array}{llllllll}
&v_1=v,\quad &v_2=J_{z_1}v,\quad &v_3=J_{z_2}v, \quad &v_4=J_{z_3}v,
\\
&v_5=J_{z_4}v,\quad &v_6=J_{z_5}v,\quad &v_7=J_{z_4}J_{z_2}v, \quad &v_8=J_{z_4}J_{z_3}v.
\end{array}
$$
{ \small
\begin{table}[h]
\caption{Commutation relations for $\n_{2,3}$}
\centering
\begin{tabular}{| c | c | c | c | c | c | c | c | c | c |} 
\hline
 $[r, c]$  & $v_1$ & $v_2$ & $v_3$ & $v_4$ & $v_5$ & $v_6$ & $v_7$ & $v_8$
 \\
\hline
$v_1$ & $0$ & $z_1$ & $z_2$ & $z_3$ & $z_4$ & $z_5$ & $0$ & $0$\\
\hline
$v_2$ & $-z_1$ & $0$ & $0$ & $0$ & $z_5$ & $-z_4$ & $-z_3$ & $-z_2$ \\
\hline
$v_3$ & $-z_2$ & $0$ & $0$ & $z_5$ & $0$ & $-z_3$ & $z_4$ & $z_1$ \\
\hline
$v_4$ & $-z_3$ & $0$ & $-z_5$ & $0$ & $0$ & $-z_2$ & $-z_1$ & $-z_4$\\
\hline
$v_5$ & $-z_4$ & $-z_5$ & $0$ & $0$ & $0$ & $-z_1$ & $z_2$ & $z_3$\\
\hline
$v_6$ & $-z_5$ & $z_4$  & $z_3$ & $z_2$ & $z_1$ & $0$ & $0$ & $0$\\
\hline
$v_7$ & $0$ & $z_3$  & $-z_4$ & $z_1$ & $-z_2$ & $0$ & $0$ & $z_5$\\
\hline
$v_8$ & $0$ & $z_2$  & $-z_1$ & $z_4$ & $-z_3$ & $0$ & $-z_5$ & $0$\\
\hline
\end{tabular}\label{23}
\end{table} }

%%%%%%%%%%%%%%%%%%%%%%%%%%%%%%%%%%%%%%%%%%%%%%%%

\subsection{ $H$-type Lie algebra $\n_{1,4}$} 

%%%%%%%%%%%%%%%%%%%%%%%%%%%%%%%%%%%%%%%%%%%%%%%%

The minimal admissible module $\vv$ is $8$-dimensional. Basis of $\mathbb R^{1,4}$ is $z_1,\ldots,z_5$ and $J_{z_1}^2=-\Id_{\vv}$, $J_{z_i}^2=\Id_{\vv}$, $i=2,3,4,5$. We choose an initial vector $v\in\vv$, $\la v,v\ra=1$ such that  
$$
P_1v=J_{z_2}J_{z_3}J_{z_4}J_{z_5}v=v,\quad P_2v=J_{z_1}J_{z_2}J_{z_3}v=v.
$$
Then the basis of $\vv$ is the following
$$
\begin{array}{llllllll}
&v_1=v,\quad &v_2=J_{z_1}v,\quad &v_3=J_{z_2}v, \quad &v_4=J_{z_3}v,
\\
&v_5=J_{z_4}v,\quad &v_6=J_{z_5}v,\quad &v_7=J_{z_4}J_{z_2}v, \quad &v_8=J_{z_4}J_{z_3}v.
\end{array}
$$
{ \small
\begin{table}[h]
\caption{Commutation relations for $\n_{1,4}$}
\centering
\begin{tabular}{| c | c | c | c | c | c | c | c | c | c |} 
\hline
 $[r, c]$  & $v_1$ & $v_2$ & $v_3$ & $v_4$ & $v_5$ & $v_6$ & $v_7$ & $v_8$
 \\
\hline
$v_1$ & $0$ & $z_1$ & $z_2$ & $z_3$ & $z_4$ & $z_5$ & $0$ & $0$\\
\hline
$v_2$ & $-z_1$ & $0$ & $z_3$ & $-z_2$ & $-z_5$ & $z_4$ & $0$ & $0$ \\
\hline
$v_3$ & $-z_2$ & $-z_3$ & $0$ & $-z_1$ & $0$ & $0$ & $-z_4$ & $z_5$ \\
\hline
$v_4$ & $-z_3$ & $z_2$ & $z_1$ & $0$ & $0$ & $0$ & $-z_5$ & $-z_4$\\
\hline
$v_5$ & $-z_4$ & $z_5$ & $0$ & $0$ & $0$ & $z_1$ & $z_2$ & $z_3$\\
\hline
$v_6$ & $-z_5$ & $-z_4$  & $0$ & $0$ & $-z_1$ & $0$ & $z_3$ & $-z_2$\\
\hline
$v_7$ & $0$ & $0$  & $z_4$ & $z_5$ & $-z_2$ & $-z_3$ & $0$ & $-z_1$\\
\hline
$v_8$ & $0$ & $0$  & $-z_5$ & $z_4$ & $-z_3$ & $z_2$ & $z_1$ & $0$\\
\hline
\end{tabular}\label{14}
\end{table} }

%%%%%%%%%%%%%%%%%%%%%%%%%%%%%%%%%%%%%%%%%%%%%%%%

\subsection{ $H$-type Lie algebra $\n_{0,5}$} 

%%%%%%%%%%%%%%%%%%%%%%%%%%%%%%%%%%%%%%%%%%%%%%%%

The minimal admissible module $\vv$ is $16$-dimensional. Basis of $\mathbb R^{0,5}$ is $z_1,\ldots,z_5$ and $J_{z_i}^2=\Id_{\vv}$, $i=1,\ldots,5$. We choose an initial vector $v\in\vv$, $\la v,v\ra=1$ in order to satisfy 
$
P_1v=J_{z_2}J_{z_3}J_{z_4}J_{z_5}v=v$.
Then the basis of $\vv$ is the following
\begin{equation*} \label{eq:basis05}
\begin{array}{lllllllllll}
& v_1=v,\quad & v_2=J_{z_1}J_{z_2}v,\quad & v_3=J_{z_1}J_{z_3}v,\quad  & v_4=J_{z_1}J_{z_4}v,
\\ 
& v_5=J_{z_1}J_{z_5}v,& v_6=J_{z_2}J_{z_5}v, & v_7=J_{z_3}J_{z_5}v, & v_8=J_{z_4}J_{z_5}v,
\\
&v_9=J_{z_1}v,\quad &v_{10}=J_{z_2}v,\quad &v_{11}=J_{z_3}v,\quad &v_{12}=J_{z_4}v,
\\ 
&v_{13}=J_{z_5}v,\ &v_{14}=J_{z_1}J_{z_2}J_{z_5}v,& v_{15}=J_{z_1}J_{z_3}J_{z_5}v,& v_{16}=J_{z_1}J_{z_4}J_{z_5}v.
\end{array}
\end{equation*}

{\tiny
\begin{table}[h]
\center\caption{Commutation relations for $\n_{0,5}$}
\begin{tabular}{| c | c | c | c | c | c | c | c | c |c | c | c | c | c | c | c | c |} 
\hline
 $[r, c]$  & $v_1$ & $v_2$ & $v_3$ & $v_4$ & $v_{5}$ & $v_{6}$ & $v_{7}$ & $v_{8}$ & $v_9$ & $v_{10}$ & $v_{11}$ & $v_{12}$ & $v_{13}$ & $v_{14}$ & $v_{15}$ & $v_{16}$  \\
\hline
$v_1$ & $0$ & $0$ & $0$ & $0$ & $0$ & $0$ & $0$ & $0$ & $z_1$ & $z_2$ & $z_3$ & $z_4$ & $z_5$ & $0$ & $0$ & $0$ \\ 
\hline
$v_2$ & $0$ & $0$ & $0$ & $0$ & $0$ & $0$ & $0$ & $0$ & $-z_2$ & $z_1$ & $0$ & $0$ & $0$ & $z_5$ & $z_4$ & $-z_3$ \\
\hline
$v_3$ & $0$ & $0$ & $0$ & $0$ & $0$ & $0$ & $0$ & $0$ & $-z_3$ & $0$ & $z_1$ & $0$ & $0$ & $-z_4$ & $z_5$ & $z_2$ \\
\hline
$v_4$ & $0$ & $0$ & $0$ & $0$ & $0$ & $0$ & $0$ & $0$ & $-z_4$ & $0$ & $0$ & $z_1$ & $0$ & $z_3$ & $-z_2$ & $z_5$ \\
\hline
$v_{5}$ & $0$ & $0$ & $0$ & $0$ & $0$ & $0$ & $0$ & $0$ & $-z_5$ & $0$ & $0$ & $0$ & $z_1$ & $-z_2$ & $-z_3$ & $-z_4$ \\
\hline
$v_{6}$ & $0$ & $0$ & $0$ & $0$ & $0$ & $0$ & $0$ & $0$ & $0$ & $-z_5$ & $z_4$ & $-z_3$ & $z_2$ & $z_1$ & $0$ & $0$ \\
\hline
$v_{7}$ & $0$ & $0$ & $0$ & $0$ & $0$ & $0$ & $0$ & $0$ & $0$ & $-z_4$ & $-z_5$ & $z_2$ & $z_3$ & $0$ & $z_1$ & $0$ \\
\hline
$v_{8}$ & $0$ & $0$ & $0$ & $0$ & $0$ & $0$ & $0$ & $0$ & $0$ & $z_3$ & $-z_2$ & $-z_5$ & $z_4$ & $0$ & $0$ & $z_1$ \\
\hline
$v_9$ & $-z_1$ & $z_2$ & $z_3$ & $z_4$ & $z_5$ & $0$ & $0$ & $0$ & $0$ & $0$ & $0$ & $0$ & $0$ & $0$ & $0$ & $0$\\
\hline
$v_{10}$ & $-z_2$ & $-z_1$ & $0$ & $0$ & $0$ & $z_5$ & $z_4$ & $-z_3$ & $0$ & $0$ & $0$ & $0$ & $0$ & $0$ & $0$ & $0$ \\
\hline
$v_{11}$ & $-z_3$ & $0$ & $-z_1$ & $0$ & $0$ & $-z_4$ & $z_5$ & $z_2$ & $0$ & $0$ & $0$ & $0$ & $0$ & $0$ & $0$ & $0$ \\
\hline
$v_{12}$ & $-z_4$ & $0$ & $0$ & $-z_1$ & $0$ & $z_3$ & $-z_2$ & $z_5$ & $0$ & $0$ & $0$ & $0$ & $0$ & $0$ & $0$ & $0$ \\
\hline
$v_{13}$ & $-z_5$ & $0$ & $0$ & $0$ & $-z_1$ & $-z_2$ & $-z_3$ & $-z_4$ & $0$ & $0$ & $0$ & $0$ & $0$ & $0$ & $0$ & $0$ \\
\hline
$v_{14}$ & $0$ & $-z_5$ & $z_4$ & $-z_3$ & $z_2$ & $-z_1$ & $0$ & $0$ & $0$ & $0$ & $0$ & $0$ & $0$ & $0$ & $0$ & $0$ \\
\hline
$v_{15}$ & $0$ & $-z_4$ & $-z_5$ & $z_2$ & $z_3$ & $0$ & $-z_1$ & $0$ & $0$ & $0$ & $0$ & $0$ & $0$ & $0$ & $0$ & $0$ \\
\hline
$v_{16}$ & $0$ & $z_3$ & $-z_2$ & $-z_5$ & $z_4$ & $0$ & $0$ & $-z_1$ & $0$ & $0$ & $0$ & $0$ & $0$ & $0$ & $0$ & $0$ \\
\hline
\end{tabular}
\label{Cl05}
%\end{sidewaystable} 
\end{table}
}

%%%%%%%%%%%%%%%%%%%%%%%%%%%%%%%%%%%%%%%%%%%%%%%%

\section{Bases and structure constants for pseudo $H$-type Lie algebras with $r+s=6$}

%%%%%%%%%%%%%%%%%%%%%%%%%%%%%%%%%%%%%%%%%%%%%%%%

%%%%%%%%%%%%%%%%%%%%%%%%%%%%%%%%%%%%%%%%%%%%%%%%

\subsection{ $H$-type Lie algebra $\n_{6,0}$} 

%%%%%%%%%%%%%%%%%%%%%%%%%%%%%%%%%%%%%%%%%%%%%%%%

the minimal admissible module $\vv$ is $8$-dimensional. Basis of $\mathbb R^{6,0}$ is $z_1,\ldots,z_6$ and $J_{z_i}^2=-\Id_{\vv}$, $i=1,\ldots,6$. We choose an initial vector $v\in\vv$, $\la v,v\ra=1$, such that
$$
P_1v=J_{z_1}J_{z_2}J_{z_3}J_{z_4}v=v,\quad P_2v=J_{z_1}J_{z_2}J_{z_5}J_{z_6}v=v,\quad P_3v=J_{z_1}J_{z_4}J_{z_5}v=v.
$$
Then the basis of $\vv$ is the following
$$
\begin{array}{llllllll}
&v_1=v,\quad &v_2=J_{z_1}v,\quad &v_3=J_{z_2}v, \quad &v_4=J_{z_3}v,
\\
&v_5=J_{z_4}v,\quad &v_6=J_{z_5}v,\quad &v_7=J_{z_6}v, \quad &v_8=J_{z_1}J_{z_2}v.
\end{array}
$$
Useful relations: $J_{z_1}J_{z_3}J_{z_6}v=v$, $-J_{z_2}J_{z_3}J_{z_5}v=v$, $J_{z_2}J_{z_4}J_{z_6}v=v$.
{ \small
\begin{table}[h]
\caption{Commutation relations for $\n_{6,0}$}
\centering
\begin{tabular}{| c | c | c | c | c | c | c | c | c | c |} 
\hline
 $[r, c]$  & $v_1$ & $v_2$ & $v_3$ & $v_4$ & $v_5$ & $v_6$ & $v_7$ & $v_8$
 \\
\hline
$v_1$ & $0$ & $z_1$ & $z_2$ & $z_3$ & $z_4$ & $z_5$ & $z_6$ & $0$\\
\hline
$v_2$ & $-z_1$ & $0$ & $0$ & $-z_6$ & $-z_5$ & $z_4$ & $z_3$ & $-z_2$ \\
\hline
$v_3$ & $-z_2$ & $0$ & $0$ & $z_5$ & $-z_6$ & $-z_3$ & $z_4$ & $z_1$ \\
\hline
$v_4$ & $-z_3$ & $z_6$ & $-z_5$ & $0$ & $0$ & $z_2$ & $-z_1$ & $z_4$\\
\hline
$v_5$ & $-z_4$ & $z_5$ & $z_6$ & $0$ & $0$ & $-z_1$ & $-z_2$ & $-z_3$\\
\hline
$v_6$ & $-z_5$ & $-z_4$ & $z_3$ & $-z_2$ & $z_1$ & $0$ & $0$ & $z_6$\\
\hline
$v_7$ & $-z_6$ & $-z_3$ & $-z_4$ & $z_1$ & $z_2$ & $0$ & $0$ & $-z_5$\\
\hline
$v_8$ & $0$ & $z_2$ & $-z_1$ & $-z_4$ & $z_3$ & $-z_6$ & $z_5$ & $0$\\
\hline
\end{tabular}\label{60}
\end{table} }

%%%%%%%%%%%%%%%%%%%%%%%%%%%%%%%%%%%%%%%%%%%%%%%%

\subsection{ $H$-type Lie algebra $\n_{5,1}$} 

%%%%%%%%%%%%%%%%%%%%%%%%%%%%%%%%%%%%%%%%%%%%%%%%

The minimal admissible module $\vv$ is $16$-dimensional. Basis of $\mathbb R^{5,1}$ is $z_1,\ldots,z_6$ and $J_{z_i}^2=-\Id_{\vv}$, $i=1,\ldots,5$, $J_{z_6}^2=\Id_{\vv}$. We choose an initial vector $v\in\vv$, $\la v,v\ra=1$, satisfying 
$$
P_1v=J_{z_1}J_{z_2}J_{z_3}J_{z_4}v=v,\quad P_2v=J_{z_1}J_{z_2}J_{z_5}v=v.
$$
Then the basis of $\vv$ for $\n_{5,1}$ is the following 
$$
\begin{array}{lllllllllllllll}
&v_1=v,\quad &v_2=J_{z_1}v,\quad &v_3=J_{z_2}v, \quad &v_4=J_{z_3}v,
\\
&v_5=J_{z_4}v,\quad &v_6=J_{z_5}v,\quad &v_7=J_{z_1}J_{z_3}v, \quad &v_8=J_{z_1}J_{z_4}v,
\\
&v_9=J_{z_6}v,\quad &v_{10}=J_{z_1}J_{z_6}v,\quad &v_{11}=J_{z_2}J_{z_6}v, \quad &v_{12}=J_{z_3}J_{z_6}v,
\\
&v_{13}=J_{z_4}J_{z_6}v,\quad &v_{14}=J_{z_5}J_{z_6}v,\quad &v_{15}=J_{z_1}J_{z_3}J_{z_6}v, \quad &v_{16}=J_{z_1}J_{z_4}J_{z_6}v.
\end{array}
$$
Useful relation: $-J_{z_1}J_{z_2}v=J_{z_5}v$.
{ \tiny
\begin{table}[h]
\caption{Commutation relations for $\n_{5,1}$}
\centering
\begin{tabular}{| c | c | c | c | c | c | c | c | c | c | c | c | c | c | c | c | c | c | c | c |} 
\hline
 $[r, c]$  & $v_1$ & $v_2$ & $v_3$ & $v_4$ & $v_5$ & $v_6$ & $v_7$ & $v_8$ & $v_9$ & $v_{10}$ & $v_{11}$ & $v_{12}$ & $v_{13}$ & $v_{14}$ & $v_{15}$ & $v_{16}$
 \\
\hline
$v_1$ & $0$ & $z_1$ & $z_2$ & $z_3$ & $z_4$ & $z_5$ & $0$ & $0$ & $z_6$ & $0$ & $0$ & $0$ & $0$ & $0$ & $0$ & $0$\\
\hline
$v_2$ & $-z_1$ & $0$ & $-z_5$ & $0$ & $0$ & $z_2$ & $-z_3$ & $-z_4$ & $0$ & $-z_6$ & $0$ & $0$ & $0$ & $0$ & $0$ & $0$\\
\hline
$v_3$ & $-z_2$ & $z_5$ & $0$ & $0$ & $0$ & $-z_1$ & $-z_4$ & $z_3$ & $0$ & $0$ & $-z_6$ & $0$ & $0$ & $0$ & $0$ & $0$\\
\hline
$v_4$ & $-z_3$ & $0$ & $0$ & $0$ & $z_5$ & $-z_4$ & $z_1$ & $-z_2$ & $0$ & $0$ & $0$ & $-z_6$ & $0$ & $0$ & $0$ & $0$\\
\hline
$v_5$ & $-z_4$ & $0$ & $0$ & $-z_5$ & $0$ & $z_3$ & $z_2$ & $z_1$ & $0$ & $0$ & $0$ & $0$ & $-z_6$ & $0$ & $0$ & $0$\\
\hline
$v_6$ & $-z_5$ & $-z_2$  & $z_1$ & $z_4$ & $-z_3$ & $0$ & $0$ & $0$ & $0$ & $0$ & $0$ & $0$ & $0$ & $-z_6$ & $0$ & $0$\\
\hline
$v_7$ & $0$ & $z_3$  & $z_4$ & $-z_1$ & $-z_2$ & $0$ & $0$ & $-z_5$ & $0$ & $0$ & $0$ & $0$ & $0$ & $0$ & $z_6$ & $0$\\
\hline
$v_8$ & $0$ & $z_4$  & $-z_3$ & $z_2$ & $-z_1$ & $0$ & $z_5$ & $0$ & $0$ & $0$ & $0$ & $0$ & $0$ & $0$ & $0$ & $z_6$\\
\hline
$v_9$ & $-z_6$ & $0$  & $0$ & $0$ & $0$ & $0$ & $0$ & $0$ & $0$ & $-z_1$  & $-z_2$ & $-z_3$ & $-z_4$ & $-z_5$ & $0$ & $0$\\
\hline
$v_{10}$ & $0$ & $z_6$ & $0$ & $0$ & $0$ & $0$ & $0$ & $0$ & $z_1$ & $0$  & $-z_5$ & $0$ & $0$ & $z_2$ & $z_3$ & $z_4$\\
\hline
$v_{11}$ & $0$ & $0$  & $z_6$ & $0$ & $0$ & $0$ & $0$ & $0$ & $z_2$ & $z_5$  & $0$ & $0$ & $0$ & $-z_1$ & $z_4$ & $-z_3$\\
\hline
$v_{12}$ & $0$ & $0$  & $0$ & $z_6$ & $0$ & $0$ & $0$ & $0$ & $z_3$ & $0$  & $0$ & $0$ & $z_5$ & $-z_4$ & $-z_1$ & $z_2$ \\
\hline
$v_{13}$ & $0$ & $0$  & $0$ & $$ & $z_6$ & $0$ & $0$ & $0$ & $z_4$ & $0$  & $0$ & $-z_5$ & $0$ & $z_3$ & $-z_2$ & $-z_1$\\
\hline
$v_{14}$ & $0$ & $0$  & $0$ & $0$ & $0$ & $z_6$ & $0$ & $0$ & $z_5$ & $-z_2$  & $z_1$ & $z_4$ & $-z_3$ & $0$ & $0$ & $0$\\
\hline
$v_{15}$ & $0$ & $0$  & $0$ & $0$ & $0$ & $0$ & $-z_6$ & $0$ & $0$ & $-z_3$  & $-z_4$ & $z_1$ & $z_2$ & $0$ & $0$ & $-z_5$\\
\hline
$v_{16}$ & $0$ & $0$  & $0$ & $0$ & $0$ & $0$ & $0$ & $-z_6$ & $0$ & $-z_4$  & $z_3$ & $-z_2$ & $z_1$ & $0$ & $z_5$ & $0$\\
\hline
\end{tabular}\label{51}
\end{table} }

%%%%%%%%%%%%%%%%%%%%%%%%%%%%%%%%%%%%%%%%%%%%%%%%

\subsection{ $H$-type Lie algebra $\n_{4,2}$} 

%%%%%%%%%%%%%%%%%%%%%%%%%%%%%%%%%%%%%%%%%%%%%%%%

The minimal admissible module $\vv$ is $16$-dimensional. Basis of $\mathbb R^{4,2}$ is $z_1,\ldots,z_6$ and $J_{z_i}^2=-\Id_{\vv}$, $i=1,2,3,4$, $J_{z_i}^2=\Id_{\vv}$, $i=5,6$. We choose an initial vector $v\in\vv$, $\la v,v\ra=1$ such that  
$$
P_1v=J_{z_1}J_{z_2}J_{z_3}J_{z_4}v=v,\quad P_2v=J_{z_1}J_{z_2}J_{z_5}J_{z_6}v=v.
$$
Then the basis of $\vv$ is the following 
$$
\begin{array}{lllllllllllllll}
&v_1=v,\quad &v_2=J_{z_1}v,\quad &v_3=J_{z_2}v, \quad &v_4=J_{z_3}v,
\\
&v_5=J_{z_4}v,\quad &v_6=J_{z_1}J_{z_2}v,\quad &v_7=J_{z_1}J_{z_3}v, \quad &v_8=J_{z_1}J_{z_4}v,
\\
&v_9=J_{z_5}v,\quad &v_{10}=J_{z_6}v,\quad &v_{11}=J_{z_1}J_{z_5}v, \quad &v_{12}=J_{z_1}J_{z_6}v,
\\
&v_{13}=J_{z_3}J_{z_5}v,\quad &v_{14}=J_{z_3}J_{z_6}v,\quad &v_{15}=J_{z_1}J_{z_3}J_{z_5}v, \quad &v_{16}=J_{z_2}J_{z_3}J_{z_5}v.
\end{array}
$$
Useful relation: $-J_{z_3}J_{z_4}J_{z_5}J_{z_6}v=v$.
{ \tiny
\begin{table}[h]
\caption{Commutation relations for $\n_{4,2}$}
\centering
\begin{tabular}{| c | c | c | c | c | c | c | c | c | c | c | c | c | c | c | c | c | c | c | c |} 
\hline
 $[r, c]$  & $v_1$ & $v_2$ & $v_3$ & $v_4$ & $v_5$ & $v_6$ & $v_7$ & $v_8$ & $v_9$ & $v_{10}$ & $v_{11}$ & $v_{12}$ & $v_{13}$ & $v_{14}$ & $v_{15}$ & $v_{16}$
 \\
\hline
$v_1$ & $0$ & $z_1$ & $z_2$ & $z_3$ & $z_4$ & $0$ & $0$ & $0$ & $z_5$ & $z_6$ & $0$ & $0$ & $0$ & $0$ & $0$ & $0$\\
\hline
$v_2$ & $-z_1$ & $0$ & $0$ & $0$ & $0$ & $-z_2$ & $-z_3$ & $-z_4$ & $0$ & $0$ & $-z_5$ & $-z_6$ & $0$ & $0$ & $0$ & $0$\\
\hline
$v_3$ & $-z_2$ & $0$ & $0$ & $0$ & $0$ & $z_1$ & $-z_4$ & $z_3$ & $0$ & $0$ & $z_6$ & $-z_5$ & $0$ & $0$ & $0$ & $0$\\
\hline
$v_4$ & $-z_3$ & $0$ & $0$ & $0$ & $0$ & $z_4$ & $z_1$ & $-z_2$ & $0$ & $0$ & $0$ & $0$ & $-z_5$ & $-z_6$ & $0$ & $0$\\
\hline
$v_5$ & $-z_4$ & $0$ & $0$ & $0$ & $0$ & $-z_3$ & $z_2$ & $z_1$ & $0$ & $0$ & $0$ & $0$ & $-z_6$ & $z_5$ & $0$ & $0$\\
\hline
$v_6$ & $0$ & $z_2$  & $-z_1$ & $-z_4$ & $z_3$ & $0$ & $0$ & $0$ & $z_6$ & $-z_5$ & $0$ & $0$ & $0$ & $0$ & $0$ & $0$\\
\hline
$v_7$ & $0$ & $z_3$  & $z_4$ & $-z_1$ & $-z_2$ & $0$ & $0$ & $0$ & $0$ & $0$ & $0$ & $0$ & $0$ & $0$ & $z_5$ & $z_6$\\
\hline
$v_8$ & $0$ & $z_4$  & $-z_3$ & $z_2$ & $-z_1$ & $0$ & $0$ & $0$ & $0$ & $0$ & $0$ & $0$ & $0$ & $0$ & $z_6$ & $-z_5$\\
\hline
$v_9$ & $-z_5$ & $0$  & $0$ & $0$ & $0$ & $-z_6$ & $0$ & $0$ & $0$ & $0$  & $-z_1$ & $-z_2$ & $-z_3$ & $z_4$ & $0$ & $0$\\
\hline
$v_{10}$ & $-z_6$ & $0$ & $0$ & $0$ & $0$ & $z_5$ & $0$ & $0$ & $0$ & $0$  & $z_2$ & $-z_1$ & $-z_4$ & $-z_3$ & $0$ & $0$\\
\hline
$v_{11}$ & $0$ & $z_5$  & $-z_6$ & $0$ & $0$ & $0$ & $0$ & $0$ & $z_1$ & $-z_2$  & $0$ & $0$ & $0$ & $0$ & $z_3$ & $-z_4$\\
\hline
$v_{12}$ & $0$ & $z_6$  & $z_5$ & $0$ & $0$ & $0$ & $0$ & $0$ & $z_2$ & $z_1$  & $0$ & $0$ & $0$ & $0$ & $z_4$ & $z_3$ \\
\hline
$v_{13}$ & $0$ & $0$  & $0$ & $z_5$ & $z_6$ & $0$ & $0$ & $0$ & $z_3$ & $z_4$  & $0$ & $0$ & $0$ & $0$ & $-z_1$ & $-z_2$\\
\hline
$v_{14}$ & $0$ & $0$  & $0$ & $z_6$ & $-z_5$ & $0$ & $0$ & $0$ & $-z_4$ & $z_3$  & $0$ & $0$ & $0$ & $0$ & $z_2$ & $-z_1$\\
\hline
$v_{15}$ & $0$ & $0$  & $0$ & $0$ & $0$ & $0$ & $-z_5$ & $-z_6$ & $0$ & $0$  & $-z_3$ & $-z_4$ & $z_1$ & $-z_2$ & $0$ & $0$\\
\hline
$v_{16}$ & $0$ & $0$  & $0$ & $0$ & $0$ & $0$ & $-z_6$ & $z_5$ & $0$ & $0$  & $z_4$ & $-z_3$ & $z_2$ & $z_1$ & $0$ & $0$\\
\hline
\end{tabular}\label{42}
\end{table} }

%%%%%%%%%%%%%%%%%%%%%%%%%%%%%%%%%%%%%%%%%%%%%%%%
%%%%%%%%%%%%%%%%%%%%%%%%%%%%%%%%%%%%%%%%%%%%%%%%

\subsection{ $H$-type Lie algebra $\n_{3,3}$} 

%%%%%%%%%%%%%%%%%%%%%%%%%%%%%%%%%%%%%%%%%%%%%%%%

The minimal admissible module $\vv$ is $8$-dimensional. Basis of $\mathbb R^{3,3}$ is $z_1,\ldots,z_6$ and $J_{z_i}^2=-\Id_{\vv}$, $i=1,2,3$, $J_{z_i}^2=\Id_{\vv}$, $i=4,5,6$. We choose an initial vector $v\in\vv$, $\la v,v\ra=1$, satisfying 
$$
P_1v=J_{z_1}J_{z_2}J_{z_4}J_{z_5}v=v,\quad P_2v=J_{z_2}J_{z_3}J_{z_5}J_{z_6}v=v,\quad P_3v=J_{z_1}J_{z_2}J_{z_3}v=v
$$
Then the basis of $\vv$ is the following
$$
\begin{array}{llllllll}
&v_1=v,\quad &v_2=J_{z_1}v,\quad &v_3=J_{z_2}v, \quad &v_4=J_{z_3}v
\\
&v_5=J_{z_4}v,\quad &v_6=J_{z_5}v,\quad &v_7=J_{z_6}v, \quad &v_8=J_{z_1}J_{z_4}v.
\end{array}
$$
Useful relations: 
$$-J_{z_1}J_{z_5}J_{z_6}v=v,\quad -J_{z_3}J_{z_4}J_{z_5}v=v,\quad -J_{z_2}J_{z_4}J_{z_6}v=v,\quad -J_{z_1}J_{z_3}J_{z_4}J_{z_6}v=v.
$$
{ \small
\begin{table}[h]
\caption{Commutation relations for $\n_{3,3}$}
\centering
\begin{tabular}{| c | c | c | c | c | c | c | c | c | c |} 
\hline
 $[r, c]$  & $v_1$ & $v_2$ & $v_3$ & $v_4$ & $v_5$ & $v_6$ & $v_7$ & $v_8$
 \\
\hline
$v_1$ & $0$ & $z_1$ & $z_2$ & $z_3$ & $z_4$ & $z_5$ & $z_6$ & $0$\\
\hline
$v_2$ & $-z_1$ & $0$ & $-z_3$ & $z_2$ & $0$ & $-z_6$ & $z_5$ & $-z_4$ \\
\hline
$v_3$ & $-z_2$ & $z_3$ & $0$ & $-z_1$ & $-z_6$ & $0$ & $z_4$ & $z_5$ \\
\hline
$v_4$ & $-z_3$ & $-z_2$ & $z_1$ & $0$ & $-z_5$ & $z_4$ & $0$ & $-z_6$\\
\hline
$v_5$ & $-z_4$ & $0$ & $z_6$ & $z_5$ & $0$ & $z_3$ & $z_2$ & $-z_1$\\
\hline
$v_6$ & $-z_5$ & $z_6$  & $0$ & $-z_4$ & $-z_3$ & $0$ & $z_1$ & $z_2$\\
\hline
$v_7$ & $-z_6$ & $-z_5$  & $-z_4$ & $0$ & $-z_2$ & $-z_1$ & $0$ & $-z_3$\\
\hline
$v_8$ & $0$ & $z_4$  & $-z_5$ & $z_6$ & $z_1$ & $-z_2$ & $z_3$ & $0$\\
\hline
\end{tabular}\label{33}
\end{table} }

%%%%%%%%%%%%%%%%%%%%%%%%%%%%%%%%%%%%%%%%%%%%%%%%

\subsection{ $H$-type Lie algebra $\n_{2,4}$} 

%%%%%%%%%%%%%%%%%%%%%%%%%%%%%%%%%%%%%%%%%%%%%%%%

The minimal admissible module $\vv$ is $8$-dimensional. Basis of $\mathbb R^{2,4}$ is $z_1,\ldots,z_6$ and $J_{z_i}^2=-\Id_{\vv}$, $i=1,2$, $J_{z_i}^2=\Id_{\vv}$, $i=3,4,5,6$. We choose an initial vector $v\in\vv$, $\la v,v\ra=1$, such that  
$$
P_1v=J_{z_1}J_{z_2}J_{z_3}J_{z_4}v=v,\quad P_2v=J_{z_1}J_{z_2}J_{z_5}J_{z_6}v=v,\quad P_3v=J_{z_1}J_{z_3}J_{z_5}v=v.
$$
Then the basis of $\vv$ for $\n_{2,4}$ is the following
$$
\begin{array}{llllllll}
&v_1=v,\quad &v_2=J_{z_1}v,\quad &v_3=J_{z_2}v, \quad &v_4=J_{z_1}J_{z_2}v,
\\
&v_5=J_{z_3}v,\quad &v_6=J_{z_4}v,\quad &v_7=J_{z_5}v, \quad &v_8=J_{z_6}v.
\end{array}
$$
Useful relations $-J_{z_2}J_{z_4}J_{z_5}v=v$, $-J_{z_1}J_{z_4}J_{z_6}v=v$, $-J_{z_2}J_{z_3}J_{z_6}v=v$.
{ \small
\begin{table}[h]
\caption{Commutation relations for $\n_{2,4}$}
\centering
\begin{tabular}{| c | c | c | c | c | c | c | c | c | c |} 
\hline
 $[r, c]$  & $v_1$ & $v_2$ & $v_3$ & $v_4$ & $v_5$ & $v_6$ & $v_7$ & $v_8$
 \\
\hline
$v_1$ & $0$ & $z_1$ & $z_2$ & $0$ & $z_3$ & $z_4$ & $z_5$ & $z_6$\\
\hline
$v_2$ & $-z_1$ & $0$ & $0$ & $-z_2$ & $z_5$ & $-z_6$ & $-z_3$ & $z_4$ \\
\hline
$v_3$ & $-z_2$ & $0$ & $0$ & $z_1$ & $-z_6$ & $-z_5$ & $z_4$ & $z_3$ \\
\hline
$v_4$ & $0$ & $z_2$ & $-z_1$ & $0$ & $z_4$ & $-z_3$ & $z_6$ & $-z_5$\\
\hline
$v_5$ & $-z_3$ & $-z_5$ & $z_6$ & $-z_4$ & $0$ & $0$ & $-z_1$ & $z_2$\\
\hline
$v_6$ & $-z_4$ & $z_6$  & $z_5$ & $z_3$ & $0$ & $0$ & $z_2$ & $z_1$\\
\hline
$v_7$ & $-z_5$ & $z_3$  & $-z_4$ & $-z_6$ & $z_1$ & $-z_2$ & $0$ & $0$\\
\hline
$v_8$ & $-z_6$ & $-z_4$  & $-z_3$ & $z_5$ & $-z_2$ & $-z_1$ & $0$ & $0$\\
\hline
\end{tabular}\label{24}
\end{table} }

%%%%%%%%%%%%%%%%%%%%%%%%%%%%%%%%%%%%%%%%%%%%%%%%

\subsection{ $H$-type Lie algebra $\n_{1,5}$} 

%%%%%%%%%%%%%%%%%%%%%%%%%%%%%%%%%%%%%%%%%%%%%%%%

The minimal admissible module $\vv$ is $16$-dimensional. Basis of $\mathbb R^{1,5}$ is $z_1,\ldots,z_6$ and $J_{z_1}^2=-\Id_{\vv}$, $J_{z_i}^2=\Id_{\vv}$, $i=2,\ldots,6$. We choose an initial vector $v\in\vv$, $\la v,v\ra=1$ satisfying 
$$
P_1v=J_{z_2}J_{z_3}J_{z_4}J_{z_5}v=v,\quad P_2v=J_{z_1}J_{z_2}J_{z_3}v=v. 
$$
Then the basis of $\vv$ is the following
\begin{equation*} \label{eq:basis15}
\begin{array}{lllllllllll}
& v_1=v,\quad & v_2=J_{z_1}v,\quad & v_3=J_{z_2}J_{z_6}v,\quad  & v_4=J_{z_3}J_{z_6}v,
\\ 
& v_5=J_{z_4}J_{z_6}v,& v_6=J_{z_5}J_{z_6}v, & v_7=J_{z_2}J_{z_4}v, & v_8=J_{z_2}J_{z_5}v,
\\
&v_9=J_{z_6}v,\quad &v_{10}=J_{z_1}J_{z_6}v,\quad &v_{11}=J_{z_2}v,\quad &v_{12}=J_{z_3}v,
\\ 
&v_{13}=J_{z_4}v,\ &v_{14}=J_{z_5}v,& v_{15}=J_{z_2}J_{z_4}J_{z_6}v,& v_{16}=J_{z_2}J_{z_5}J_{z_6}v.
\end{array}
\end{equation*}
 Useful relation: $-J_{z_1}J_{z_4}J_{z_5}v=v$.
{\tiny
\begin{table}[h]
\center\caption{Commutation relations for $\n_{1,5}$}
\begin{tabular}{| c | c | c | c | c | c | c | c | c |c | c | c | c | c | c | c | c |} 
\hline
 $[r, c]$  & $v_1$ & $v_2$ & $v_3$ & $v_4$ & $v_{5}$ & $v_{6}$ & $v_{7}$ & $v_{8}$ & $v_9$ & $v_{10}$ & $v_{11}$ & $v_{12}$ & $v_{13}$ & $v_{14}$ & $v_{15}$ & $v_{16}$  \\
\hline
$v_1$ & $0$ & $z_1$ & $0$ & $0$ & $0$ & $0$ & $0$ & $0$ & $z_6$ & $0$ & $z_2$ & $z_3$ & $z_4$ & $z_5$ & $0$ & $0$ \\ 
\hline
$v_2$ & $-z_1$ & $0$ & $0$ & $0$ & $0$ & $0$ & $0$ & $0$ & $0$ & $-z_6$ & $z_3$ & $-z_2$ & $-z_5$ & $z_4$ & $0$ & $0$ \\
\hline
$v_3$ & $0$ & $0$ & $0$ & $-z_1$ & $0$ & $0$ & $0$ & $0$ & $z_2$ & $-z_3$ & $-z_6$ & $0$ & $0$ & $0$ & $-z_4$ & $-z_5$ \\
\hline
$v_4$ & $0$ & $0$ & $z_1$ & $0$ & $0$ & $0$ & $0$ & $0$ & $z_3$ & $z_2$ & $0$ & $-z_6$ & $0$ & $0$ & $-z_5$ & $z_4$ \\
\hline
$v_{5}$ & $0$ & $0$ & $0$ & $0$ & $0$ & $z_1$ & $0$ & $0$ & $z_4$ & $z_5$ & $0$ & $0$ & $-z_6$ & $0$ & $z_2$ & $-z_3$ \\
\hline
$v_{6}$ & $0$ & $0$ & $0$ & $0$ & $-z_1$ & $0$ & $0$ & $0$ & $z_5$ & $-z_4$ & $0$ & $0$ & $0$ & $-z_6$ & $z_3$ & $z_2$ \\
\hline
$v_{7}$ & $0$ & $0$ & $0$ & $0$ & $0$ & $0$ & $0$ & $z_1$ & $0$ & $0$ & $-z_4$ & $-z_5$ & $z_2$ & $z_3$ & $z_6$ & $0$ \\
\hline
$v_{8}$ & $0$ & $0$ & $0$ & $0$ & $0$ & $0$ & $-z_1$ & $0$ & $0$ & $0$ & $-z_5$ & $z_4$ & $-z_3$ & $z_2$ & $0$ & $z_6$ \\
\hline
$v_9$ & $-z_6$ & $0$ & $-z_2$ & $-z_3$ & $-z_4$ & $-z_5$ & $0$ & $0$ & $0$ & $-z_1$ & $0$ & $0$ & $0$ & $0$ & $0$ & $0$\\
\hline
$v_{10}$ & $0$ & $z_6$ & $z_3$ & $-z_2$ & $-z_5$ & $z_4$ & $0$ & $0$ & $z_1$ & $0$ & $0$ & $0$ & $0$ & $0$ & $0$ & $0$ \\
\hline
$v_{11}$ & $-z_2$ & $-z_3$ & $z_6$ & $0$ & $0$ & $0$ & $z_4$ & $z_5$ & $0$ & $0$ & $0$ & $-z_1$ & $0$ & $0$ & $0$ & $0$ \\
\hline
$v_{12}$ & $-z_3$ & $z_2$ & $0$ & $z_6$ & $0$ & $0$ & $z_5$ & $-z_4$ & $0$ & $0$ & $z_1$ & $0$ & $0$ & $0$ & $0$ & $0$ \\
\hline
$v_{13}$ & $-z_4$ & $z_5$ & $0$ & $0$ & $z_6$ & $0$ & $-z_2$ & $z_3$ & $0$ & $0$ & $0$ & $0$ & $0$ & $z_1$ & $0$ & $0$ \\
\hline
$v_{14}$ & $-z_5$ & $-z_4$ & $0$ & $0$ & $0$ & $z_6$ & $-z_3$ & $-z_2$ & $0$ & $0$ & $0$ & $0$ & $-z_1$ & $0$ & $0$ & $0$ \\
\hline
$v_{15}$ & $0$ & $0$ & $z_4$ & $z_5$ & $-z_2$ & $-z_3$ & $-z_6$ & $0$ & $0$ & $0$ & $0$ & $0$ & $0$ & $0$ & $0$ & $z_1$ \\
\hline
$v_{16}$ & $0$ & $0$ & $z_5$ & $-z_4$ & $z_3$ & $-z_2$ & $0$ & $-z_6$ & $0$ & $0$ & $0$ & $0$ & $0$ & $0$ & $-z_1$ & $0$ \\
\hline
\end{tabular}
\label{15}
%\end{sidewaystable} 
\end{table}
}

%%%%%%%%%%%%%%%%%%%%%%%%%%%%%%%%%%%%%%%%%%%%%%%%

\subsection{ $H$-type Lie algebra $\n_{0,6}$} 

%%%%%%%%%%%%%%%%%%%%%%%%%%%%%%%%%%%%%%%%%%%%%%%%

The minimal admissible module $\vv$ is $16$-dimensional. Basis of $\mathbb R^{0,6}$ is $z_1,\ldots,z_6$ and $J_{z_i}^2=\Id_{\vv}$, $i=1,\ldots, 6$. We choose an initial vector $v\in\vv$, $\la v,v\ra=1$ satisfying 
$$
P_1v=J_{z_1}J_{z_2}J_{z_3}J_{z_4}v=v,\quad P_2v=J_{z_1}J_{z_2}J_{z_5}J_{z_6}v.
$$
Then the basis of $\vv$ is the following
\begin{equation*} \label{eq:basis05}
\begin{array}{lllllllllll}
& v_1=v,\quad & v_2=J_{z_1}J_{z_2}v,\quad & v_3=J_{z_1}J_{z_3}v,\quad  & v_4=J_{z_1}J_{z_4}v,
\\ 
& v_5=J_{z_1}J_{z_5}v,& v_6=J_{z_1}J_{z_6}v, & v_7=J_{z_3}J_{z_5}v, & v_8=J_{z_3}J_{z_6}v,
\\
&v_9=J_{z_1}v,\quad &v_{10}=J_{z_2}v,\quad &v_{11}=J_{z_3}v,\quad &v_{12}=J_{z_4}v,
\\ 
&v_{13}=J_{z_5}v,\ &v_{14}=J_{z_6}v,& v_{15}=J_{z_1}J_{z_3}J_{z_5}v,& v_{16}=J_{z_1}J_{z_3}J_{z_6}v.
\end{array}
\end{equation*}
 Useful relation:  $-J_{z_3}J_{z_4}J_{z_5}J_{z_6}v=v$.
{\tiny
\begin{table}[h]
\center\caption{Commutation relations for $\n_{0,6}$}
\begin{tabular}{| c | c | c | c | c | c | c | c | c |c | c | c | c | c | c | c | c |} 
\hline
 $[r, c]$  & $v_1$ & $v_2$ & $v_3$ & $v_4$ & $v_{5}$ & $v_{6}$ & $v_{7}$ & $v_{8}$ & $v_9$ & $v_{10}$ & $v_{11}$ & $v_{12}$ & $v_{13}$ & $v_{14}$ & $v_{15}$ & $v_{16}$  \\
\hline
$v_1$ & $0$ & $0$ & $0$ & $0$ & $0$ & $0$ & $0$ & $0$ & $z_1$ & $z_2$ & $z_3$ & $z_4$ & $z_5$ & $z_6$ & $0$ & $0$ \\ 
\hline
$v_2$ & $0$ & $0$ & $0$ & $0$ & $0$ & $0$ & $0$ & $0$ & $-z_2$ & $z_1$ & $z_4$ & $-z_3$ & $z_6$ & $-z_5$ & $0$ & $0$ \\
\hline
$v_3$ & $0$ & $0$ & $0$ & $0$ & $0$ & $0$ & $0$ & $0$ & $-z_3$ & $-z_4$ & $z_1$ & $z_2$ & $0$ & $0$ & $z_5$ & $z_6$ \\
\hline
$v_4$ & $0$ & $0$ & $0$ & $0$ & $0$ & $0$ & $0$ & $0$ & $-z_4$ & $z_3$ & $-z_2$ & $z_1$ & $0$ & $0$ & $-z_6$ & $z_5$ \\
\hline
$v_{5}$ & $0$ & $0$ & $0$ & $0$ & $0$ & $0$ & $0$ & $0$ & $-z_5$ & $-z_6$ & $0$ & $0$ & $z_1$ & $z_2$ & $-z_3$ & $-z_4$ \\
\hline
$v_{6}$ & $0$ & $0$ & $0$ & $0$ & $0$ & $0$ & $0$ & $0$ & $-z_6$ & $z_5$ & $0$ & $0$ & $-z_2$ & $z_1$ & $z_4$ & $-z_3$ \\
\hline
$v_{7}$ & $0$ & $0$ & $0$ & $0$ & $0$ & $0$ & $0$ & $0$ & $0$ & $0$ & $-z_5$ & $z_6$ & $z_3$ & $-z_4$ & $z_1$ & $-z_2$ \\
\hline
$v_{8}$ & $0$ & $0$ & $0$ & $0$ & $0$ & $0$ & $0$ & $0$ & $0$ & $0$ & $-z_6$ & $-z_5$ & $z_4$ & $z_3$ & $z_2$ & $z_1$ \\
\hline
$v_9$ & $-z_1$ & $z_2$ & $z_3$ & $z_4$ & $z_5$ & $z_6$ & $0$ & $0$ & $0$ & $0$ & $0$ & $0$ & $0$ & $0$ & $0$ & $0$\\
\hline
$v_{10}$ & $-z_2$ & $-z_1$ & $z_4$ & $-z_3$ & $z_6$ & $-z_5$ & $0$ & $0$ & $0$ & $0$ & $0$ & $0$ & $0$ & $0$ & $0$ & $0$ \\
\hline
$v_{11}$ & $-z_3$ & $-z_4$ & $-z_1$ & $z_2$ & $0$ & $0$ & $z_5$ & $z_6$ & $0$ & $0$ & $0$ & $0$ & $0$ & $0$ & $0$ & $0$ \\
\hline
$v_{12}$ & $-z_4$ & $z_3$ & $-z_2$ & $-z_1$ & $0$ & $0$ & $-z_6$ & $z_5$ & $0$ & $0$ & $0$ & $0$ & $0$ & $0$ & $0$ & $0$ \\
\hline
$v_{13}$ & $-z_5$ & $-z_6$ & $0$ & $0$ & $-z_1$ & $z_2$ & $-z_3$ & $-z_4$ & $0$ & $0$ & $0$ & $0$ & $0$ & $0$ & $0$ & $0$ \\
\hline
$v_{14}$ & $-z_6$ & $z_5$ & $0$ & $0$ & $-z_2$ & $-z_1$ & $z_4$ & $-z_3$ & $0$ & $0$ & $0$ & $0$ & $0$ & $0$ & $0$ & $0$ \\
\hline
$v_{15}$ & $0$ & $0$ & $-z_5$ & $z_6$ & $z_3$ & $-z_4$ & $-z_1$ & $-z_2$ & $0$ & $0$ & $0$ & $0$ & $0$ & $0$ & $0$ & $0$ \\
\hline
$v_{16}$ & $0$ & $0$ & $-z_6$ & $-z_5$ & $z_4$ & $z_3$ & $z_2$ & $-z_1$ & $0$ & $0$ & $0$ & $0$ & $0$ & $0$ & $0$ & $0$ \\
\hline
\end{tabular}
\label{Cl06}
%\end{sidewaystable} 
\end{table}
}

%%%%%%%%%%%%%%%%%%%%%%%%%%%%%%%%%%%%%%%%%%%%%%%%

\section{Bases and structure constants for pseudo $H$-type Lie algebras with $r+s=7$} 

%%%%%%%%%%%%%%%%%%%%%%%%%%%%%%%%%%%%%%%%%%%%%%%%

%%%%%%%%%%%%%%%%%%%%%%%%%%%%%%%%%%%%%%%%%%%%%%%%

\subsection{ $H$-type Lie algebra $\n_{7,0}$} 

%%%%%%%%%%%%%%%%%%%%%%%%%%%%%%%%%%%%%%%%%%%%%%%%

The minimal admissible module $\vv$ is $8$-dimensional. Basis of $\mathbb R^{7,0}$ is $z_1,\ldots,z_7$ and $J_{z_i}^2=-\Id_{\vv}$, $i=1,\ldots,7$. We choose an initial vector $v\in\vv$, $\la v,v\ra=1$, such that 
$$
\begin{array}{lllll}
&P_1v=J_{z_1}J_{z_2}J_{z_3}J_{z_4}v=v,\quad &P_2v=J_{z_1}J_{z_2}J_{z_5}J_{z_6}v=v,
\\
&P_3v=J_{z_1}J_{z_3}J_{z_5}J_{z_7}v=v,\quad &P_4v=J_{z_5}J_{z_6}J_{z_7}v=v.
\end{array} 
$$
The basis of $\vv$ is the following 
$$
\begin{array}{llllll}
& v_1=v, \quad &v_2=J_{z_1}v, \quad &v_3=J_{z_2}v, \quad &v_4=J_{z_3}v, 
\\ 
& v_5=J_{z_4}v, \quad &v_6=J_{z_5}v, \quad &v_7=J_{z_6}v, \quad &v_8=J_{z_7}v,.
\end{array}
$$
Useful relations: 
$$P_2P_4v=-J_{z_1}J_{z_2}J_{z_7}v=v, \quad P_3P_4v=J_{z_1}J_{z_3}J_{z_6}v=v,\quad P_1P_2P_3P_4v=J_{z_1}J_{z_4}J_{z_5}v=v,
$$
$$
P_2P_3P_4v=-J_{z_2}J_{z_3}J_{z_5}v=v,\quad P_1P_3P_4v=J_{z_2}J_{z_4}J_{z_6}v=v,\quad P_1P_2P_4v=J_{z_3}J_{z_4}J_{z_7}v=v.
$$

{ \small
\begin{table}[h]
\caption{Commutation relations for $\n_{7,0}$}
\centering
\begin{tabular}{| c | c | c | c | c | c | c | c | c | c |} 
\hline
 $[r, c]$  & $v_1$ & $v_2$ & $v_3$ & $v_4$ & $v_5$ & $v_6$ & $v_7$ & $v_8$
 \\
\hline
$v_1$ & $0$ & $z_1$ & $z_2$ & $z_3$ & $z_4$ & $z_5$ & $z_6$ & $z_7$\\
\hline
$v_2$ & $-z_1$ & $0$ & $z_7$ & $-z_6$ & $-z_5$ & $z_4$ & $z_3$ & $-z_2$ \\
\hline
$v_3$ & $-z_2$ & $-z_7$ & $0$ & $z_5$ & $-z_6$ & $-z_3$ & $z_4$ & $z_1$ \\
\hline
$v_4$ & $-z_3$ & $z_6$ & $-z_5$ & $0$ & $-z_7$ & $z_2$ & $-z_1$ & $z_4$\\
\hline
$v_5$ & $-z_4$ & $z_5$ & $z_6$ & $z_7$ & $0$ & $-z_1$ & $-z_2$ & $-z_3$\\
\hline
$v_6$ & $-z_5$ & $-z_4$  & $z_3$ & $-z_2$ & $z_1$ & $0$ & $-z_7$ & $z_6$\\
\hline
$v_7$ & $-z_6$ & $-z_3$  & $-z_4$ & $z_1$ & $z_2$ & $z_7$ & $0$ & $-z_5$\\
\hline
$v_8$ & $-z_7$ & $z_2$  & $-z_1$ & $-z_4$ & $z_3$ & $-z_6$ & $z_5$ & $0$\\
\hline
\end{tabular}\label{70}
\end{table} }

%%%%%%%%%%%%%%%%%%%%%%%%%%%%%%%%%%%%%%%%%%%%%%%%

\subsection{ $H$-type Lie algebra $\n_{3,4}$} 

%%%%%%%%%%%%%%%%%%%%%%%%%%%%%%%%%%%%%%%%%%%%%%%%

The minimal admissible module $\vv$ is $8$-dimensional. Basis of $\mathbb R^{3,4}$ is $z_1,\ldots,z_7$ with $J_{z_i}^2=-\Id_{\vv}$, $i=1,\ldots,3$ and $J_{z_i}^2=\Id_{\vv}$, $i=4,\ldots,7$. We choose an initial vector $v\in\vv$, $\la v,v\ra=1$, such that 
$$
\begin{array}{lllll}
&P_1v=J_{z_1}J_{z_2}J_{z_4}J_{z_5}v=v,\quad &P_2v=J_{z_1}J_{z_2}J_{z_6}J_{z_7}v=v,
\\
&P_3v=J_{z_1}J_{z_3}J_{z_5}J_{z_7}v=v,\quad &P_4v=J_{z_1}J_{z_2}J_{z_3}v=v.
\end{array} 
$$
The basis of $\vv$ is the following 
$$
\begin{array}{llllll}
& v_1=v, \quad &v_2=J_{z_1}v, \quad &v_3=J_{z_2}v, \quad &v_4=J_{z_3}v, 
\\ 
& v_5=J_{z_4}v, \quad &v_6=J_{z_5}v, \quad &v_7=J_{z_6}v, \quad &v_8=J_{z_7}v,.
\end{array}
$$
Useful relations: 
$$P_1P_3P_4v=-J_{z_1}J_{z_4}J_{z_7}v=v, \quad P_2P_3P_4v=-J_{z_1}J_{z_5}J_{z_6}v=v,\quad P_1P_2P_3P_4v=-J_{z_2}J_{z_4}J_{z_6}v=v,
$$
$$
P_3P_4v=J_{z_2}J_{z_5}J_{z_7}v=v,\quad P_1P_4v=-J_{z_3}J_{z_4}J_{z_5}v=v,\quad P_2P_4v=-J_{z_3}J_{z_6}J_{z_7}v=v.
$$

{ \small
\begin{table}[h]
\caption{Commutation relations for $\n_{3,4}$}
\centering
\begin{tabular}{| c | c | c | c | c | c | c | c | c | c |} 
\hline
 $[r, c]$  & $v_1$ & $v_2$ & $v_3$ & $v_4$ & $v_5$ & $v_6$ & $v_7$ & $v_8$
 \\
\hline
$v_1$ & $0$ & $z_1$ & $z_2$ & $z_3$ & $z_4$ & $z_5$ & $z_6$ & $z_7$\\
\hline
$v_2$ & $-z_1$ & $0$ & $-z_3$ & $z_2$ & $-z_7$ & $-z_6$ & $z_5$ & $z_4$ \\
\hline
$v_3$ & $-z_2$ & $z_3$ & $0$ & $-z_1$ & $-z_6$ & $z_7$ & $z_4$ & $-z_5$ \\
\hline
$v_4$ & $-z_3$ & $-z_2$ & $z_1$ & $0$ & $-z_5$ & $z_4$ & $-z_7$ & $z_6$\\
\hline
$v_5$ & $-z_4$ & $z_7$ & $z_6$ & $z_5$ & $0$ & $z_3$ & $z_2$ & $z_1$\\
\hline
$v_6$ & $-z_5$ & $z_6$  & $-z_7$ & $-z_4$ & $-z_3$ & $0$ & $z_1$ & $-z_2$\\
\hline
$v_7$ & $-z_6$ & $-z_5$  & $-z_4$ & $z_7$ & $-z_2$ & $-z_1$ & $0$ & $z_3$\\
\hline
$v_8$ & $-z_7$ & $-z_4$  & $z_5$ & $-z_6$ & $-z_1$ & $z_2$ & $-z_3$ & $0$\\
\hline
\end{tabular}\label{34}
\end{table} }

%%%%%%%%%%%%%%%%%%%%%%%%%%%%%%%%%%%%%%%%%%%%%%%%

\subsection{ $H$-type Lie algebra $\n_{0,7}$} 

%%%%%%%%%%%%%%%%%%%%%%%%%%%%%%%%%%%%%%%%%%%%%%%%

The minimal admissible module $\vv$ is $16$-dimensional. Basis of $\mathbb R^{0,7}$ is $z_1,\ldots,z_7$ with $J_{z_i}^2=\Id_{\vv}$, $i=1,\ldots,7$. A result from~\cite[Theorem11]{FM2} states that the pseudo $H$-type Lie algebra $\n_{0,7}$ is isomorphic to the pseudo $H$-type Lie algebra $\n_{7,0}(\vv^{+}\oplus \vv^-)$ which minimal admissible module is the direct sum of two minimal admissible modules for $\n_{7,0}$. Here $\vv^{+}$ and $\vv^{-}$ are two non-equivalent irreducible modules of the Clifford algebra $\Cl_{7,0}$.
The direct sum $\vv^{+}\oplus \vv^-$ is orthogonal and therefore the table of structural constants will have block diagonal form. To calculate the first block we consider four mutually commuting involutions action on In the first block we

Let $z_{1},\ldots,z_7, ~\text{are positive orthonormal generators}$ of
$C\ell_{7,0}$ and
let mutually commuting positive involutions $P_{i}$ 
be as above for the admissible module $V_{min}^{+}\cong\mathbb{R}^{8,0}$ 
of $C\ell_{7,0}$, that is
\begin{align*}
&P_{1}=J_{Z_1}J_{Z_2}J_{Z_3}J_{Z_4},~P_{2}=J_{Z_1}J_{Z_2}J_{Z_5}J_{Z_6},\\
&P_{3}=J_{Z_1}J_{Z_3}J_{Z_5}J_{Z_7},~P_{4}=J_{Z_5}J_{Z_6}J_{Z_7},
\end{align*}
and for the admissible module $V_{min}^{-}$ 
\begin{align*}
&Q_{1}=\tilde{J}_{Z_1}\tilde{J}_{Z_2}\tilde{J}_{Z_3}\tilde{J}_{Z_4},
~Q_{2}=\tilde{J}_{Z_1}\tilde{J}_{Z_2}\tilde{J}_{Z_5}\tilde{J}_{Z_6},\\
&Q_{3}=\tilde{J}_{Z_1}\tilde{J}_{Z_3}\tilde{J}_{Z_5}\tilde{J}_{Z_7},
~Q_{4}=\tilde{J}_{Z_5}\tilde{J}_{Z_6}\tilde{J}_{Z_7}.
\end{align*}

For the admissible module $V_{min}^{+}$ of $C\ell_{7,0}$ 
we fix a vector $v\in V_{min}^{+}$
with the properties that
\[
P_{i}(v)=v, i=1,2,3,4~\text{and}~<v,v>_{8,0}=1
\]
and for another admissible module $V_{min}^{-}\cong\mathbb{R}^{8,0}$ of $C\ell_{7,0}$
we take a vector $w\in V_{min}^{-}$ with the properties that
\[
P_{i}(w)=w, i=1,2,3~text{and}~P_{4}(w)=-w,~\text{and}~<w,w>_{8,0}=1.
\]

We choose positive orthonormal basis $\{X_{i}\}_{i=0}^{15}$ of $V_{min}^{+}\oplus V_{min}^{-}$ with
\begin{align*}
&X_{0}=v,X_{1}=J_{z_1}(v),X_{2}=J_{z_2}(v),X_{3}=J_{z_3}(v),X_{4}=J_{z_4}(v),\\
&X_{5}=J_{z_5}(v),X_{6}=J_{z_6}(v),X_{7}=J_{z_7}(v),\\
&X_{8}= w, X_{9}=\tilde{J}_{z_1}(w),
X_{10}=\tilde{J}_{z_2}(w),X_{11}=\tilde{J}_{z_3}(w),X_{12}=\tilde{J}_{z_4}(w),\\
&X_{13}=\tilde{J}_{z_5}(w),X_{14}=\tilde{J}_{z_6}(w),X_{15}=\tilde{J}_{z_7}(w).
\end{align*}

Then it holds always $[X_i, X_{8+j}]=0$ for $i,j=0,7$.

%%%%%%%%%%%%%%%%%%%%%%%%%%%%%%%%%%%%%%%%%%%%%%%%

\section{Bases and structure constants for pseudo $H$-type Lie algebras with $r+s=8$} 

%%%%%%%%%%%%%%%%%%%%%%%%%%%%%%%%%%%%%%%%%%%%%%%%

%%%%%%%%%%%%%%%%%%%%%%%%%%%%%%%%%%%%%%%%%%%%%%%%

\subsection{ $H$-type Lie algebra $\n_{8,0}$} 

%%%%%%%%%%%%%%%%%%%%%%%%%%%%%%%%%%%%%%%%%%%%%%%%

The minimal admissible module $\vv$ is $16$-dimensional. Basis of $\mathbb R^{8,0}$ is $z_1,\ldots,z_8$ and $J_{z_i}^2=-\Id_{\vv}$, $i=1,\ldots,8$. We choose an initial vector $v\in\vv$, $\la v,v\ra=1$, such that 
$$
\begin{array}{lllll}
&P_1v=J_{z_1}J_{z_2}J_{z_3}J_{z_4}v=v,\quad &P_2v=J_{z_1}J_{z_2}J_{z_5}J_{z_6}v=v,
\\
&P_3v=J_{z_2}J_{z_3}J_{z_5}J_{z_7}v=v,\quad &P_4v=J_{z_1}J_{z_2}J_{z_7}J_{z_8}v=v.
\end{array} 
$$
The basis of $\vv$ is the following 
$$
\begin{array}{llllll}
& v_1=v, \quad &v_2=J_{z_1}J_{z_2}v, \quad &v_3=J_{z_1}J_{z_3}v, \quad &v_4=J_{z_1}J_{z_4}v, 
\\ 
& v_5=J_{z_1}J_{z_5}v, \quad &v_6=J_{z_1}J_{z_6}v, \quad &v_7=J_{z_1}J_{z_7}v, \quad &v_8=J_{z_1}J_{z_8}v,
\\
& v_9=J_{z_1}v, \quad &v_{10}=J_{z_2}v, \quad &v_{11}=J_{z_3}v, \quad &v_{12}=J_{z_4}v, 
\\ 
& v_{13}=J_{z_5}v, \quad &v_{14}=J_{z_6}v, \quad &v_{15}=J_{z_7}v, \quad &v_{16}=J_{z_8}v.
\end{array}
$$
Useful relations: 
$$P_3P_4v=-J_{z_1}J_{z_3}J_{z_5}J_{z_8}v=v,\quad P_2P_3v=-J_{z_1}J_{z_3}J_{z_6}J_{z_7}v=v,
$$
$$
P_1P_3v=-J_{z_1}J_{z_4}J_{z_5}J_{z_7}v=v,\quad P_1P_2P_3P_4v=J_{z_1}J_{z_4}J_{z_6}J_{z_8}v=v.
$$
The pseudo $H$-type algebras $\n_{8,0}$ and $\n_{8,0}$ are isomorphic
\begin{table}[h!]
{\tiny 
\center\caption{Commutation relations for $\n_{8,0}$ and $\n_{0,8}$}
\begin{tabular}{|c|c| c | c | c | c | c | c | c |c | c | c | c | c | c | c | c |}
\hline
$[r,c]$&$v_1$&$v_2$ & $v_3$ & $v_4$ & $v_5$ & $v_6$ & $v_7$ & $v_8$ & $v_9$ & $v_{10}$ & $v_{11}$ & $v_{12}$ & $v_{13}$ & $v_{14}$ & $v_{15}$ & $v_{16}$  \\
\hline
$v_1$ & $0$ & $0$ & $0$ & $0$ & $0$ & $0$ & $0$ & $0$ & $z_1$ & $z_2$ & $z_3$ & $z_4$ & $z_5$ & $z_6$ & $z_7$ & $z_8$ \\ 
\hline
$v_2$ & $0$ & $0$ & $0$ & $0$ & $0$ & $0$ & $0$ & $0$ & $z_2$ & $-z_1$ & $-z_4$ & $z_3$ & $-z_6$ & $z_5$ & $-z_8$ & $z_7$ \\
\hline
$v_3$ & $0$ & $0$ & $0$ & $0$ & $0$ & $0$ & $0$ & $0$ & $z_3$ & $z_4$ & $-z_1$ & $-z_2$ & $z_8$ & $z_7$ & $-z_6$ & $-z_5$ \\
\hline
$v_4$ & $0$ & $0$ & $0$ & $0$ & $0$ & $0$ & $0$ & $0$ & $z_4$ & $-z_3$ & $z_2$ & $-z_1$ & $z_7$ & $-z_8$ & $-z_5$ & $z_6$ \\
\hline
$v_5$ & $0$ & $0$ & $0$ & $0$ & $0$ & $0$ & $0$ & $0$ & $z_5$ & $z_6$ & $-z_8$ & $-z_7$ & $-z_1$ & $-z_2$ & $z_4$ & $z_3$ \\
\hline
$v_6$ & $0$ & $0$ & $0$ & $0$ & $0$ & $0$ & $0$ & $0$ & $z_6$ & $-z_5$ & $-z_7$ & $z_8$ & $z_2$ & $-z_1$ & $z_3$ & $-z_4$ \\
\hline
$v_7$ & $0$ & $0$ & $0$ & $0$ & $0$ & $0$ & $0$ & $0$ & $z_7$ & $z_8$ & $z_6$ & $z_5$ & $-z_4$ & $-z_3$ & $-z_1$ & $-z_2$ \\
\hline
$v_8$ & $0$ & $0$ & $0$ & $0$ & $0$ & $0$ & $0$ & $0$ & $z_8$ & $-z_7$ & $z_5$ & $-z_6$ & $-z_3$ & $z_4$ & $z_2$ & $-z_1$ \\
\hline
$v_9$ & $-z_1$ & $-z_2$ & $-z_3$ & $-z_4$ & $-z_5$ & $-z_6$ & $-z_7$ & $-z_8$ & $0$ & $0$ & $0$ & $0$ & $0$ & $0$ & $0$ & $0$\\
\hline
$v_{10}$ & $-z_2$ & $z_1$ & $-z_4$ & $z_3$ & $-z_6$ & $z_5$ & $-z_8$ & $z_7$ & $0$ & $0$ & $0$ & $0$ & $0$ & $0$ & $0$ & $0$ \\
\hline
$v_{11}$ & $-z_3$ & $z_4$ & $z_1$ & $-z_2$ & $z_8$ & $z_7$ & $-z_6$ & $-z_5$ & $0$ & $0$ & $0$ & $0$ & $0$ & $0$ & $0$ & $0$ \\
\hline
$v_{12}$ & $-z_4$ & $-z_3$ & $z_2$ & $z_1$ & $z_7$ & $-z_8$ & $-z_5$ & $z_6$ & $0$ & $0$ & $0$ & $0$ & $0$ & $0$ & $0$ & $0$ \\
\hline
$v_{13}$ & $-z_5$ & $z_6$ & $-z_8$ & $-z_7$ & $z_1$ & $-z_2$ & $z_4$ & $z_3$ & $0$ & $0$ & $0$ & $0$ & $0$ & $0$ & $0$ & $0$ \\
\hline
$v_{14}$ & $-z_6$ & $-z_5$ & $-z_7$ & $z_8$ & $z_2$ & $z_1$ & $z_3$ & $-z_4$ & $0$ & $0$ & $0$ & $0$ & $0$ & $0$ & $0$ & $0$ \\
\hline
$v_{15}$ & $-z_7$ & $z_8$ & $z_6$ & $z_5$ & $-z_4$ & $-z_3$ & $z_1$ & $-z_2$ & $0$ & $0$ & $0$ & $0$ & $0$ & $0$ & $0$ & $0$ \\
\hline
$v_{16}$ & $-z_8$ & $-z_7$ & $z_5$ & $-z_6$ & $-z_3$ & $z_4$ & $z_2$ & $z_1$ & $0$ & $0$ & $0$ & $0$ & $0$ & $0$ & $0$ & $0$ \\
\hline 
\end{tabular}
\label{80}
}
\end{table}

%%%%%%%%%%%%%%%%%%%%%%%%%%%%%%%%%%%%%%%%%%%%%%%%

\subsection{ $H$-type Lie algebra $\n_{7,1}$} 

%%%%%%%%%%%%%%%%%%%%%%%%%%%%%%%%%%%%%%%%%%%%%%%%

The minimal admissible module $\vv$ is $16$-dimensional. Basis of $\mathbb R^{7,1}$ is $z_1,\ldots,z_8$ with $J_{z_i}^2=-\Id_{\vv}$, $i=1,\ldots,7$ and $J_{z_8}^2=\Id_{\vv}$. We choose an initial vector $v\in\vv$, $\la v,v\ra=1$, such that 
$$
\begin{array}{lllll}
&P_1v=J_{z_1}J_{z_2}J_{z_3}J_{z_4}v=v,\quad &P_2v=J_{z_1}J_{z_2}J_{z_5}J_{z_6}v=v,
\\
&P_3v=J_{z_1}J_{z_3}J_{z_5}J_{z_7}v=v,\quad &P_4v=J_{z_5}J_{z_6}J_{z_7}=v.
\end{array} 
$$
The basis of $\vv$ is the following 
$$
\begin{array}{llllll}
& v_1=v, \quad &v_2=J_{z_1}v, \quad &v_3=J_{z_2}v, \quad &v_4=J_{z_3}v, 
\\ 
& v_5=J_{z_4}v, \quad &v_6=J_{z_5}v, \quad &v_7=J_{z_6}v, \quad &v_8=J_{z_7}v,
\\
& v_9=J_{z_8}v, \quad &v_{10}=J_{z_8}J_{z_1}v, \quad &v_{11}=J_{z_8}J_{z_2}v, \quad &v_{12}=J_{z_8}J_{z_3}v, 
\\ 
& v_{13}=J_{z_8}J_{z_4}v, \quad &v_{14}=J_{z_8}J_{z_5}v, \quad &v_{15}=J_{z_8}J_{z_6}v, \quad &v_{16}=J_{z_8}J_{z_7}v.
\end{array}
$$
Useful relations: 
$$P_2P_4v=-J_{z_1}J_{z_2}J_{z_7}v=v, \quad P_3P_4v=J_{z_1}J_{z_3}J_{z_6}v=v,\quad P_1P_2P_3P_4v=J_{z_1}J_{z_4}J_{z_5}v=v,
$$
$$
P_2P_3P_4v=-J_{z_2}J_{z_3}J_{z_5}v=v,\quad P_1P_3P_4v=J_{z_2}J_{z_4}J_{z_6}v=v,\quad P_1P_2P_4v=J_{z_3}J_{z_4}J_{z_7}v=v.
$$
\begin{table}[h!]
{\tiny 
\center\caption{Commutation relations for $\n_{7,1}$}
\begin{tabular}{|c|c| c | c | c | c | c | c | c |c | c | c | c | c | c | c | c |}
\hline
$[r,c]$&$v_1$&$v_2$ & $v_3$ & $v_4$ & $v_5$ & $v_6$ & $v_7$ & $v_8$ & $v_9$ & $v_{10}$ & $v_{11}$ & $v_{12}$ & $v_{13}$ & $v_{14}$ & $v_{15}$ & $v_{16}$  \\
\hline
$v_1$ & $0$ & $z_1$ & $z_2$ & $z_3$ & $z_4$ & $z_5$ & $z_6$ & $z_7$ & $z_8$ & $0$ & $0$ & $0$ & $0$ & $0$ & $0$ & $0$  \\ 
\hline
$v_2$ & $-z_1$ & $0$ & $z_7$ & $-z_6$ & $-z_5$ & $z_4$ & $z_3$ & $-z_2$ & $0$ & $z_8$ & $0$ & $0$ & $0$ & $0$ & $0$ & $0$  \\
\hline
$v_3$ & $-z_2$ & $-z_7$ & $0$ & $z_5$ & $-z_6$ & $-z_3$ & $z_4$ & $z_1$ & $0$ & $0$ & $z_8$& $0$ & $0$ & $0$ & $0$ & $0$  \\
\hline
$v_4$ & $-z_3$ & $z_6$ & $-z_5$ & $0$ & $-z_7$ & $z_2$ & $-z_1$ & $z_4$ & $0$ & $0$ & $0$ & $z_8$ & $0$ & $0$ & $0$ & $0$  \\
\hline
$v_5$ & $-z_4$ & $z_5$ & $z_6$ & $z_7$ & $0$ & $-z_1$ & $-z_2$ & $-z_3$ & $0$ & $0$ & $0$ & $0$ & $z_8$ & $0$ & $0$ & $0$  \\
\hline
$v_6$ & $-z_5$ & $-z_4$ & $z_3$ & $-z_2$ & $z_1$ & $0$ & $-z_7$ & $z_6$ & $0$ & $0$ & $0$ & $0$ & $0$& $z_8$ & $0$ & $0$  \\
\hline
$v_7$ & $-z_6$ & $z_3$ & $-z_4$ & $z_1$ & $z_2$ & $z_7$ & $0$ & $-z_5$ & $0$ & $0$ & $0$ & $0$ & $0$ & $0$ & $z_8$ & $0$ \\
\hline
$v_8$ & $-z_7$ & $z_2$ & $-z_1$ & $-z_4$ & $z_3$ & $-z_6$& $z_5$ & $0$ & $0$ & $0$ & $0$ & $0$ & $0$ & $0$  & $0$ & $z_8$ \\
\hline
$v_9$ & $-z_8$ & $0$ & $0$ & $0$ & $0$ & $0$ & $0$ & $0$ & $0$ & $z_1$ & $z_2$ & $z_3$ & $z_4$ & $z_5$ & $z_6$ & $z_7$\\
\hline
$v_{10}$ & $0$ & $-z_8$ & $0$ & $0$ & $0$ & $0$ & $0$ & $0$ & $-z_1$ & $0$ & $z_7$ & $-z_6$ & $-z_5$ & $z_4$ & $z_3$ & $-z_2$  \\
\hline
$v_{11}$ & $0$ & $0$ & $-z_8$ & $0$ & $0$ & $0$ & $0$ & $0$ & $-z_2$ & $-z_7$ & $0$ & $z_5$ & $-z_6$ & $-z_3$ & $z_4$ & $z_1$  \\
\hline
$v_{12}$ & $0$ & $0$ & $0$ & $-z_8$ & $0$ & $0$ & $0$ & $0$ & $-z_3$ & $z_6$ & $-z_5$ & $0$ & $-z_7$ & $z_2$ & $-z_1$ & $z_4$  \\
\hline
$v_{13}$ & $0$ & $0$ & $0$ & $0$ & $-z_8$ & $0$ & $0$ & $0$ & $-z_4$ & $z_5$ & $z_6$ & $z_7$ & $0$ & $-z_1$ & $-z_2$ & $-z_3$ \\
\hline
$v_{14}$ & $0$ & $0$ & $0$ & $0$ & $0$ & $-z_8$ & $0$ & $0$ & $-z_5$ & $-z_4$ & $z_3$ & $-z_2$ & $z_1$ & $0$ & $-z_7$ & $z_6$ \\
\hline
$v_{15}$ & $0$ & $0$ & $0$ & $0$ & $0$ & $0$ & $-z_8$ & $0$ & $-z_6$ & $-z_3$ & $-z_4$ & $z_1$ & $z_2$ & $z_7$ & $0$ & $-z_5$ \\
\hline
$v_{16}$ & $0$ & $0$ & $0$ & $0$ & $0$ & $0$ & $0$ & $-z_8$ & $-z_7$ & $z_2$ & $-z_1$ & $-z_4$ & $z_3$ & $-z_6$ & $z_5$ & $0$ \\
\hline 
\end{tabular}
\label{71}
}
\end{table}

%%%%%%%%%%%%%%%%%%%%%%%%%%%%%%%%%%%%%%%%%%%%%%%%

\subsection{ $H$-type Lie algebra $\n_{4,4}$} 

%%%%%%%%%%%%%%%%%%%%%%%%%%%%%%%%%%%%%%%%%%%%%%%%

The minimal admissible module is $16$-dimensional. Basis of $\mathbb R^{4,4}$ is $z_1,\ldots,z_8$ and $J_{z_i}^2=-\Id_{\vv}$, $i=1,2,3,4$, $J_{z_i}^2=\Id_{\vv}$, $i=5,6,7,8$. We choose an initial vector $v\in\vv$, $\la v,v\ra=1$, such that 
$$
\begin{array}{lllllll}
&P_1v=J_{z_1}J_{z_2}J_{z_3}J_{z_4}v=v,\quad &P_2v =J_{z_1}J_{z_2}J_{z_5}J_{z_6}v=v,
\\
&P_3v=J_{z_2}J_{z_3}J_{z_5}J_{z_7}v=v,\ &P_4v=J_{z_1}J_{z_2}J_{z_7}J_{z_8}v=v.
\end{array}
$$
Then the basis of $\vv$ is the following
$$
\begin{array}{llllllll}
&v_1=v,\quad &v_2=J_{z_1}J_{z_2}v,\quad &v_3=J_{z_1}J_{z_3}v, \quad &v_4=J_{z_1}J_{z_4}v,
\\
&v_5=J_{z_1}J_{z_5}v,\quad &v_6=J_{z_1}J_{z_6}v,\quad &v_7=J_{z_1}J_{z_7}v, \quad &v_8=J_{z_1}J_{z_8}v,
\\
&v_9=J_{z_1}v,\quad &v_{10}=J_{z_2}v,\quad &v_{11}=J_{z_3}v, \quad &v_{12}=J_{z_4}v,
\\
&v_{13}=J_{z_5}v,\quad &v_{14}=J_{z_6}v,\quad &v_{15}=J_{z_7}v, \quad &v_{16}=J_{z_8}v.
\end{array}
$$
Useful relations: 
$$P_3P_4v=J_{z_1}J_{z_3}J_{z_5}J_{z_8}v=v,\quad P_2P_3v=J_{z_1}J_{z_3}J_{z_6}J_{z_7}v=v,$$
$$ P_1P_3v=-J_{z_1}J_{z_4}J_{z_5}J_{z_7}v=v,\ P_1P_3P_3P_4v=J_{z_1}J_{z_4}J_{z_6}J_{z_8}v=v.
$$
{\tiny
\begin{table}
\center\caption{Commutation relations for $\n_{4,4}$}
\begin{tabular}{| c | c | c | c | c | c | c | c | c |c | c | c | c | c | c | c | c |} 
\hline
 $[r, c]$  & $v_1$ & $v_2$ & $v_3$ & $v_4$ & $v_{5}$ & $v_{6}$ & $v_{7}$ & $v_{8}$ & $v_9$ & $v_{10}$ & $v_{11}$ & $v_{12}$ & $v_{13}$ & $v_{14}$ & $v_{15}$ & $v_{16}$  \\
\hline
$v_1$ & $0$ & $0$ & $0$ & $0$ & $0$ & $0$ & $0$ & $0$ & $z_1$ & $z_2$ & $z_3$ & $z_4$ & $z_5$ & $z_6$ & $z_7$ & $z_8$ \\ 
\hline
$v_2$ & $0$ & $0$ & $0$ & $0$ & $0$ & $0$ & $0$ & $0$ & $z_2$ & $-z_1$ & $-z_4$ & $z_3$ & $z_6$ & $-z_5$ & $z_8$ & $-z_7$ \\
\hline
$v_3$ & $0$ & $0$ & $0$ & $0$ & $0$ & $0$ & $0$ & $0$ & $z_3$ & $z_4$ & $-z_1$ & $-z_2$ & $z_8$ & $z_7$ & $-z_6$ & $-z_5$ \\
\hline
$v_4$ & $0$ & $0$ & $0$ & $0$ & $0$ & $0$ & $0$ & $0$ & $z_4$ & $-z_3$ & $z_2$ & $-z_1$ & $-z_7$ & $z_8$ & $z_5$ & $-z_6$ \\
\hline
$v_{5}$ & $0$ & $0$ & $0$ & $0$ & $0$ & $0$ & $0$ & $0$ & $z_5$ & $-z_6$ & $-z_8$ & $z_7$ & $z_1$ & $-z_2$ & $z_4$ & $-z_3$ \\
\hline
$v_{6}$ & $0$ & $0$ & $0$ & $0$ & $0$ & $0$ & $0$ & $0$ & $z_6$ & $z_5$ & $-z_7$ & $-z_8$ & $z_2$ & $z_1$ & $-z_3$ & $-z_4$ \\
\hline
$v_{7}$ & $0$ & $0$ & $0$ & $0$ & $0$ & $0$ & $0$ & $0$ & $z_7$ & $-z_8$ & $z_6$ & $-z_5$ & $-z_4$ & $z_3$ & $z_1$ & $-z_2$ \\
\hline
$v_{8}$ & $0$ & $0$ & $0$ & $0$ & $0$ & $0$ & $0$ & $0$ & $z_8$ & $z_7$ & $z_5$ & $z_6$ & $z_3$ & $z_4$ & $z_2$ & $z_1$ \\
\hline
$v_9$ & $-z_1$ & $-z_2$ & $-z_3$ & $-z_4$ & $-z_5$ & $-z_6$ & $-z_7$ & $-z_8$ & $0$ & $0$ & $0$ & $0$ & $0$ & $0$ & $0$ & $0$\\
\hline
$v_{10}$ & $-z_2$ & $z_1$ & $-z_4$ & $z_3$ & $z_6$ & $-z_5$ & $z_8$ & $-z_7$ & $0$ & $0$ & $0$ & $0$ & $0$ & $0$ & $0$ & $0$ \\
\hline
$v_{11}$ & $-z_3$ & $z_4$ & $z_1$ & $-z_2$ & $z_8$ & $z_7$ & $-z_6$ & $-z_5$ & $0$ & $0$ & $0$ & $0$ & $0$ & $0$ & $0$ & $0$ \\
\hline
$v_{12}$ & $-z_4$ & $-z_3$ & $z_2$ & $z_1$ & $-z_7$ & $z_8$ & $z_5$ & $-z_6$ & $0$ & $0$ & $0$ & $0$ & $0$ & $0$ & $0$ & $0$ \\
\hline
$v_{13}$ & $-z_5$ & $-z_6$ & $-z_8$ & $z_7$ & $-z_1$ & $-z_2$ & $z_4$ & $-z_3$ & $0$ & $0$ & $0$ & $0$ & $0$ & $0$ & $0$ & $0$ \\
\hline
$v_{14}$ & $-z_6$ & $z_5$ & $-z_7$ & $-z_8$ & $z_2$ & $-z_1$ & $-z_3$ & $-z_4$ & $0$ & $0$ & $0$ & $0$ & $0$ & $0$ & $0$ & $0$ \\
\hline
$v_{15}$ & $-z_7$ & $-z_8$ & $z_6$ & $-z_5$ & $-z_4$ & $z_3$ & $-z_1$ & $-z_2$ & $0$ & $0$ & $0$ & $0$ & $0$ & $0$ & $0$ & $0$ \\
\hline
$v_{16}$ & $-z_8$ & $z_7$ & $z_5$ & $z_6$ & $z_3$ & $z_4$ & $z_2$ & $-z_1$ & $0$ & $0$ & $0$ & $0$ & $0$ & $0$ & $0$ & $0$ \\
\hline
\end{tabular}
\label{Cl44}
\end{table}
}

%%%%%%%%%%%%%%%%%%%%%%%%%%%%%%%%%%%%%%%%%%%%%%%%

\subsection{ $H$-type Lie algebra $\n_{3,5}$} 

%%%%%%%%%%%%%%%%%%%%%%%%%%%%%%%%%%%%%%%%%%%%%%%%

The minimal admissible module $\vv$ is $16$-dimensional. Basis of $\mathbb R^{3,5}$ is $z_1,\ldots,z_7$ with $J_{z_i}^2=-\Id_{\vv}$, $i=1,\ldots,3$ and $J_{z_i}^2=\Id_{\vv}$, $i=4,\ldots,8$. We choose an initial vector $v\in\vv$, $\la v,v\ra=1$, such that 
$$P_1v=J_{z_1}J_{z_2}J_{z_4}J_{z_5}v=v,\qquad
P_2v=J_{z_1}J_{z_2}J_{z_6}J_{z_7}v=v,
$$
$$
P_3v=J_{z_1}J_{z_3}J_{z_5}J_{z_7}v=v,\qquad P_4v=J_{z_1}J_{z_2}J_{z_3}v=v.
$$
The basis for $\vv$ is the following 
$$
\begin{array}{llllll}
& v_1=v, \quad &v_2=J_{z_1}v, \quad &v_3=J_{z_2}v, \quad &v_4=J_{z_3}v, 
\\ 
& v_5=J_{z_4}v, \quad &v_6=J_{z_5}v, \quad &v_7=J_{z_6}v, \quad &v_8=J_{z_7}v,
\\
& v_9=J_{z_8}v, \quad &v_{10}=J_{z_8}J_{z_1}v, \quad &v_{11}=J_{z_8}J_{z_2}v, \quad &v_{12}=J_{z_8}J_{z_3}v, 
\\ 
& v_{13}=J_{z_8}J_{z_4}v, \quad &v_{14}=J_{z_8}J_{z_5}v, \quad &v_{15}=J_{z_8}J_{z_6}v, \quad &v_{16}=J_{z_8}J_{z_7}v.
\end{array}
$$
Useful relations: $P_3P_4v=J_{z_2}J_{z_5}J_{z_7}v=v$,\ \  $P_1P_4v=-J_{z_3}J_{z_4}J_{z_5}v=v$, 
$$P_2P_4v=-J_{z_3}J_{z_6}J_{z_7}v=v,\quad P_1P_3P_4v=-J_{z_1}J_{z_4}J_{z_7}v=v, 
$$
$$ P_2P_3P_4v=-J_{z_1}J_{z_5}J_{z_6}v=v,\quad
 P_1P_2P_3P_4v=-J_{z_2}J_{z_4}J_{z_6}v=v.
$$
{ \tiny
\begin{table}[h]
\caption{Commutation relations for $\n_{3,5}$}
\centering
\begin{tabular}{|c|c| c | c | c | c | c | c | c |c | c | c | c | c | c | c | c |}
\hline
 $[r, c]$  & $v_1$ & $v_2$ & $v_3$ & $v_4$ & $v_5$ & $v_6$ & $v_7$ & $v_8$ & $v_9$ & $v_{10}$ & $v_{11}$ & $v_{12}$ & $v_{13}$ & $v_{14}$ & $v_{15}$ & $v_{16}$
 \\
\hline
$v_1$ & $0$ & $z_1$ & $z_2$ & $z_3$ & $z_4$ & $z_5$ & $z_6$ & $z_7$ & $z_8$ & $0$ & $0$ & $0$ & $0$ & $0$ & $0$ & $0$\\
\hline
$v_2$ & $-z_1$ & $0$ & $-z_3$ & $z_2$ & $-z_7$ & $-z_6$ & $z_5$ & $z_4$ & $0$ & $z_8$ & $0$ & $0$ & $0$ & $0$ & $0$ & $0$\\
\hline
$v_3$ & $-z_2$ & $z_3$ & $0$ & $-z_1$ & $-z_6$ & $z_7$ & $z_4$ & $-z_5$ & $0$ & $0$ & $z_8$& $0$ & $0$ & $0$ & $0$ & $0$ \\
\hline
$v_4$ & $-z_3$ & $-z_2$ & $z_1$ & $0$ & $-z_5$ & $z_4$ & $-z_7$ & $z_6$ & $0$ & $0$ & $0$ & $z_8$ & $0$ & $0$ & $0$ & $0$\\
\hline
$v_5$ & $-z_4$ & $z_7$ & $z_6$ & $z_5$ & $0$ & $z_3$ & $z_2$ & $z_1$ & $0$ & $0$ & $0$ & $0$ & $-z_8$ & $0$ & $0$ & $0$\\
\hline
$v_6$ & $-z_5$ & $z_6$  & $-z_7$ & $-z_4$ & $-z_3$ & $0$ & $z_1$ & $-z_2$ & $0$ & $0$ & $0$ & $0$ & $0$& $-z_8$ & $0$ & $0$\\
\hline
$v_7$ & $-z_6$ & $-z_5$  & $-z_4$ & $z_7$ & $-z_2$ & $-z_1$ & $0$ & $z_3$ & $0$ & $0$ & $0$ & $0$ & $0$ & $0$ & $-z_8$ & $0$ \\
\hline
$v_8$ & $-z_7$ & $-z_4$  & $z_5$ & $-z_6$ & $-z_1$ & $z_2$ & $-z_3$ & $0$ & $0$ & $0$  & $0$ & $0$ & $0$ & $0$  & $0$ & $-z_8$\\
\hline
$v_9$ & $-z_8$ & $0$ & $0$ & $0$ & $0$ & $0$ & $0$ & $0$ & $0$ & $z_1$ & $z_2$ & $z_3$ & $z_4$ & $z_5$ & $z_6$ & $z_7$\\
\hline
$v_{10}$ & $0$ & $-z_8$ & $0$ & $0$ & $0$ & $0$ & $0$ & $0$ & $-z_1$ & $0$ & $-z_3$ & $z_2$ & $-z_7$ & $-z_6$ & $z_5$ & $z_4$  \\
\hline
$v_{11}$ & $0$ & $0$ & $-z_8$ & $0$ & $0$ & $0$ & $0$ & $0$ & $-z_2$ & $z_3$ & $0$ & $-z_1$ & $-z_6$ & $z_7$ & $z_4$ & $-z_5$ \\
\hline
$v_{12}$ & $0$ & $0$ & $0$ & $-z_8$ & $0$ & $0$ & $0$ & $0$ & $-z_3$ & $-z_2$ & $z_1$ & $0$ & $-z_5$ & $z_4$ & $-z_7$ & $z_6$ \\
\hline
$v_{13}$ & $0$ & $0$ & $0$ & $0$ & $z_8$ & $0$ & $0$ & $0$ & $-z_4$ & $z_7$ & $z_6$ & $z_5$ & $0$ & $z_3$ & $z_2$ & $z_1$ \\
\hline
$v_{14}$ & $0$ & $0$ & $0$ & $0$ & $0$ & $z_8$ & $0$ & $0$ & $-z_5$ & $z_6$ & $-z_7$ & $-z_4$ & $-z_3$ & $0$ & $z_1$ & $-z_2$ \\
\hline
$v_{15}$ & $0$ & $0$ & $0$ & $0$ & $0$ & $0$ & $z_8$ & $0$ & $-z_6$ & $-z_5$ & $-z_4$ & $z_7$ & $-z_2$ & $-z_1$ & $0$ & $z_3$ \\
\hline
$v_{16}$ & $0$ & $0$ & $0$ & $0$ & $0$ & $0$ & $0$ & $z_8$ & $-z_7$ & $-z_4$ & $z_5$ & $-z_6$ & $-z_1$ & $z_2$ & $z_3$ & $0$ \\
\hline
\end{tabular}\label{35}
\end{table} }

%%%%%%%%%%%%%%%%%%%%%%%%%%%%%%%%%%%%%%%%%%%%%%%%

\section{Table of Clifford algebras}

%%%%%%%%%%%%%%%%%%%%%%%%%%%%%%%%%%%%%%%%%%%%%%%%
By the red colour we denote the Clifford algebras having minimal admissible module as a double of its irreducible.
\begin{table}[h]
\center\caption{Clifford algebras}
\begin{tabular}{|c||c|c|c|c|c|c|c|c|c|c|}
\hline
${\text{\small 9}} $&$ \vdots$&$ \vdots$&$ \vdots$&$ \vdots$&$ \vdots$&$ \vdots$&$ \vdots$&$\vdots$&$ \vdots$&$\ldots$
\\
\hline
${\text{\small 8}} $&$\mathbb R(16)$&$\mathbb C(16)$&$\mathbb H(16)
$&$\mathbb H^2(16)$&$\mathbb H(32)$&$\mathbb C(64)$&$\mathbb R(128)$&$
\mathbb R^2(128)$&$\mathbb R(256)$&$\ldots$
\\
\hline
${\text{\small 7}}$ &$ \mathbb C(8)$&$\mathbb H(8)$&$
{\color{red}\mathbb H^2(8)} $&$\mathbb H(16)
$&${\color{red}\mathbb C(32)}$&${\color{red}\mathbb R(64)}  $&$\color{red}\mathbb R^2(64)$&$\mathbb R(128) $&$\mathbb C(128)$&$\ldots$
\\
\hline
${\text{\small 6}}$ &$\mathbb H(4)$&$\mathbb H^2(4)$&$\mathbb H(8)$&$\mathbb C(16)$&${\color{red}\mathbb R(32)}
$&$\color{red}\mathbb R^2(32)$&${\color{red}\mathbb R(64)} $&$\mathbb C(64)$&$ \mathbb H(64)$&$\ldots$
\\
\hline
${\text{\small 5}} $&$\color{red}\mathbb H^2(2)$&$\mathbb H(4)$&$\mathbb C(8)$&$\mathbb R(16)$&$\color{red}\mathbb R^2(16)$&$\color{red}\mathbb R(32) $&$\color{red}\mathbb C(32)$&$\mathbb H(32)$&${\color{red}\mathbb H^2(32)}$&$\ldots$
\\
\hline
${\text{\small 4}} $&${\small\mathbb H(2)}$&${\small\mathbb C(4)}$&$\mathbb R(8)
$&$ \mathbb R^2(8)$&$\mathbb R(16)$&$\mathbb C(16)$&$\mathbb H(16)
$&$\mathbb H^2(16)$&$\mathbb H(32)$&$\ldots$
\\
\hline
${\text{\small 3}}$&${\color{red} \mathbb C(2)}$&${\color{red} \mathbb R(4)}$&$\color{red} \mathbb R^2(4)$&$ \mathbb R(8)$&$ \mathbb C(8)$&$ \mathbb H(8)$
&$\color{red} \mathbb H^2(8)$&$ \mathbb H(16)$&${\color{red}\ \mathbb C(32)}$&$\ldots$
\\
\hline
${\text{\small 2}}$&${\color{red} \mathbb R(2)}$&$
{\color{red} \mathbb R^2(2)}$&${\color{red} \mathbb R(4)}$&$  \mathbb C(4)$&$ \mathbb H(4)$&$ \mathbb H^2(4)$&$ \mathbb H(8)$&$ \mathbb C(16) $&${\color{red} \mathbb R(32)}$&$\ldots$
\\
\hline
${\text{\small 1}}$ &${\color{red} {\small\mathbb R^2}}$&${\color{red} \mathbb R(2)}$&${\color{red} \mathbb C(2)}$& $ \mathbb H(2)$&${\color{red}\mathbb H^2(2)}$&$ \mathbb H(4)$&$ \mathbb C(8)$&$ \mathbb R(16)$&${\color{red}{\mathbb R^2(16)}}$ &$\ldots$
\\
\hline
${\text{\small 0}} $&$  \mathbb R$&$  \mathbb C$&$  \mathbb H$&$  \mathbb H^2$&$  \mathbb H(2)$&$  \mathbb C(4)$&$  \mathbb R(8)$&$  \mathbb R(8)$&$ \mathbb R(16)$&$\ldots$
\\
\hline\hline
{s/r}&  {\text{\small 0}}& {\text{\small 1}}& 
{\text{\small 2}}&{\text{\small 3}} & {\text{\small 4}}& {\text{\small 5}}& {\text{\small 6}}& {\text{\small 7}}& {\text{\small 8}}&$\ldots$
\\
\hline
\end{tabular}\label{t:dim}
\end{table}

%%%%%%%%%%%%%%%%%%%%%%%%%%%%%%%%%%%%%%%%%%%%%%%%%%%%%%%%

\end{document}